\documentclass[11pt]{amsart}

\usepackage{amsmath,amssymb,amscd}
\usepackage[dvips]{graphics}

\newtheorem{theorem}{Theorem}[section]
\newtheorem{corollary}[theorem]{Corollary}

\newtheorem{lemma}[theorem]{Lemma}
\newtheorem{proposition}[theorem]{Proposition}
\newtheorem{example}[theorem]{Example}
\newtheorem{definition}[theorem]{Definition}
\newtheorem{remark}[theorem]{Remark}

\def\PP{\mathbb{P}}

\newcommand{\cE}{{\mathcal E}}

\newcommand{\cL}{{\mathcal L}}
\newcommand{\cO}{{\mathcal O}}

\newcommand{\cJ}{{\mathcal J}}
\newcommand{\cK}{{\mathcal K}}

\def\GL{\operatorname{GL}}
\def\pf{\operatorname{Pf}}
\def\ker{\operatorname{Ker}}
\def\coker{\operatorname{Coker}}

\def\div{\operatorname{div}}
\def\Id{\operatorname{Id}}
\def\pic{\operatorname{Pic}}
\def\ext{\operatorname{Ext}}
\def\rank{\operatorname{rank}}
\def\elem{\operatorname{elem}}
\def\trace{\operatorname{Trace}}

\begin{document}

\title{Elementary Transformations of Pfaffian Representations of Plane Curves}
\author{Anita Buckley}

\address{Department of Mathematics, University of Ljubljana,
Jadranska 19, 1000 Ljubljana, Slovenia  
{\rm e-mail:} {\it anita.buckley@fmf.uni-lj.si}}%

\begin{abstract}
Let $C$ be a smooth curve in $\PP^2$ given by an equation $F=0$ of degree $d$. 
In this paper we consider elementary transformations of linear pfaffian representations of $C$. 
Elementary transformations can be interpreted as actions on a rank 2 vector bundle
on $C$ with canonical determinant and no sections, which corresponds to the cokernel of a pfaffian representation. 
Every two pfaffian representations of $C$ can be bridged by a finite sequence of elementary transformations. 
Pfaffian representations and elementary transformations are constructed explicitly. For a smooth quartic, applications to Aronhold bundles and 
theta characteristics are given.
\end{abstract}

\maketitle

\section{Introduction}
\label{introdef}

Let $k$ be an algebraically closed field and $C$ a nonsingular curve defined by an irreducible polynomial $F(x_0,x_1,x_2)$ 
of degree $d$ in $\PP^2$. 
A \textit{linear pfaffian representation} of $C$ is a $2d\times 2d$ skew-symmetric matrix 
$$A= \left[\begin{array}{ccccc}
0       & L_{1\, 2} & L_{1\, 3}   & \cdots & L_{1\, 2d}\\
-L_{1\, 2} & 0        & L_{2\, 3} & \cdots & L_{2\, 2d}\\
-L_{1\, 3}  & -L_{2\, 3}& 0        &        &          \\
 \vdots   &  \vdots &          & \ddots &  \vdots         \\
 -L_{1\, 2d} & -L_{2\, 2d} &          &  \cdots  & 0
\end{array}\right]$$
with linear forms $L_{ij}=a_{ij}^0 x_0+a_{ij}^1 x_1+a_{ij}^2 x_2$
such that 
$$\pf A(x_0,x_1,x_2)=c\, F(x_0,x_1,x_2)\ \mbox{ for some }\ c\in k, c\neq 0.$$
Its cokernel is a rank 2 vector bundle on $C$. Throughout the paper we equate the notion of 
vector bundles and locally free sheaves.  

Two pfaffian representations $A$ and $A'$ are \textit{equivalent} if there exists $X\in\GL_{2d}(k)$ such that
$$A'=XAX^t.$$

There is a one to one correspondence between linear pfaffian representations (up to equivalence) of $C$ and
rank 2 vector bundles (up to isomorphism) on $C$ with certain properties. This well known result is summed up in the 
following theorem of Beauville~\cite[Corollary 2.4]{beauville}.

\begin{theorem}\label{beauvilcor}
Let $C$ be a smooth plane curve defined by a polynomial $F$ of degree $d$ and let $\cE$ be a rank 2 bundle on $C$ 
with determinant $\cO_C(d-1)$ and $H^0(C, \cE(-1))=0$. 
Then there exists a $2d\times 2d$ skew-symmetric linear matrix $A$ with $\pf A=F$ and an exact sequence 
\begin{eqnarray}\label{sxsqrk1}
0\rightarrow \bigoplus_{i=1}^{2d} \cO_{\PP^2}(-1)\stackrel{A}{\longrightarrow} \bigoplus_{i=1}^{2d} \cO_{\PP^2}
\rightarrow \cE \rightarrow 0.
\end{eqnarray}

Conversely, let $A$ be a linear skew-symmetric $2d\times 2d$ matrix with $\pf A=F$. 
Then its cokernel is a a rank 2 bundle  
with $\det \cE\cong \cO_C(d-1)$ and $H^0(C, \cE(-1))=0$.  
\end{theorem}
Using this, in~\cite{anitapfaff} all linear pfaffian representations of $C$ (up to equivalence) were found and related to 
\textit{the moduli space} $M_C(2,K_C)$
of semistable rank $2$ vector bundles on $C$ with canonical determinant. 
In particular, pfaffian representations of $C$ can be parametrised by the open set  $M_C(2,K_C)-\{\mathcal{K}\, :\, h^0(C, \mathcal{K})>0 \}$.
The properties of the moduli space were extensively studied in~\cite{newstead},~\cite{seshadri} and ~\cite{potier}; for example it
is an irreducible, normal projective variety
and for $C$ of genus $g\geq 2$ it has dimension $3(g-1)$. 

Study of pfaffian representations is strongly related to and motivated by determinantal representations. 
A \textit{linear determinantal representation} of $C$ is a $d\times d$  matrix 
$M(x_0,x_1,x_2)$ of linear forms 
such that 
$$\det M(x_0,x_1,x_2)=c\, F(x_0,x_1,x_2)\ \mbox{ for some }\ c\in k, c\neq 0.$$
Determinantal representations $M$ and $M'$ are \textit{equivalent} if there exists $X,Y\in\GL_{d}(k)$ such that
$$M'=XMY.$$
By~\cite{vinnikov2} all linear determinantal representations of $C$ (up to equivalence) can be parametrised by the open set in the Jacobian
variety 
$$\{\mbox{line bundle }\mathcal{L}:\, \deg\mathcal{L}=\frac{1}{2}d(d-3),\, h^0(C, \cL)=0 \}.$$
Determinantal representations can be seen as a special case of pfaffian representations.
Indeed, every determinantal representation $M$ 
induces a \textit{ decomposable pfaffian representation}
\begin{equation}\label{decomrep}
\left[ \begin{array}{cc}
0 & M \\
-M^t & 0
\end{array}\right].
\end{equation}
The corresponding cokernel equals $\coker M \oplus (\coker^{-1} M \otimes \cO(d-1))$ as described 
in~\cite{anitapfaff} and~\cite {dolgachev}. In the moduli space, decomposable pfaffian representations 
correspond to an open subset of the singular locus of $M_C(2,K_C)$.
Note that the equivalence relation is also well defined since
$$\left[ \begin{array}{cc}
0 & X M Y \\
-(X M Y)^t & 0
		\end{array}\right]=\left[ \begin{array}{cc}
X & 0 \\
0 & Y^t
\end{array}\right]\left[ \begin{array}{cc}
0 & M \\
-M^t & 0
\end{array}\right]\left[ \begin{array}{cc}
X^t & 0 \\
0 & Y
\end{array}\right].$$ 

Elementary transformations of determinantal representations were introduced by M. S. Livsic and Kravitsky 
in operator theory. Vinnikov et al generalised these ideas in~\cite{vinn1},~\cite{vinn2} using 
notions of vessels and Cauchy kernels. An explicit and complete description of elementary transformations of determinantal 
representations can be found  in~\cite{vinnikovElTr}.
Elementary transformations of vector bundles are due to Maruyama~\cite{maruyama}. For a modern exposition and proofs we refer to Abe~\cite{abe}.
The most traditional example of elementary transformations is that of ruled surfaces (i.e., of rank 2 bundles).
A good reference on ruled surfaces is~\cite[Chapter V. 2]{hartshorne}.

Throughout the paper we will be using the following properties of the wedge product.
It is well known that $\bigwedge^2 k^{2d}$ can be identified 
with $2d\times 2d$ skew-symmetric matrices. Let $\Omega_n$ denote the set of vectors
$\sum_{i=1}^{n} w_i\wedge w'_i$ in $\bigwedge^2 k^{2d}$, where\\ $\dim \{w_1,\ldots,w_n,w'_1,\ldots,w'_n \}=2n$. 
Elements of 
$\Omega_n$ are said to have \textit{irreducible length} $n$ since they can be written as a sum of $n$ and not 
less than $n$ pure nonzero products.
In~\cite{westwick} it is shown that 
$\Omega_n$ is isomorphic to the set of all rank $2n$ skew-symmetric matrices. 
The isomorphism equals 
\begin{equation}
\label{explrk}
\sum_{i=1}^{2d} \alpha_i e_i\ \wedge\  \sum_{j=1}^{2d} \beta_j e_j \ \  \mapsto \ \ 
\sum_{i,j=1}^{2d}(\alpha_i \beta_j-\alpha_j \beta_i)(E_{ij}-E_{ji}),
\end{equation}
where $\{e_1,\ldots,e_{2d}\}$ is the standard basis for $k^{2d}$ and $\{E_{ij}\}$ is the standard basis for $2d\times 2d$ matrices.
Note that these are the Pl\"{u}cker coordinates in Gr$(2,2d)$.

A brief outline of the paper is the following. In section~\ref{canform} we present pfaffian representations in canonical forms. This enables us to 
write an algorithm which computes all the pfaffian representations of $C$ and thus gives a description of $M_C(2,K_C)$.
For suitable choices of vectors in the cokernel bundles discussed in Section~\ref {tangents}, we define elementary transformations of pfaffian
representations in Section~\ref {beginsec}. In Section~\ref{secthree} elementary transformations of pfaffian
representations are related to elementary transformations of vector bundles. The main Theorem~\ref{cortwoparts} in Section~\ref{decsec} proves that 
any two pfaffian representations can be bridged by a finite sequence of elementary transformations of Type I and II. In other words, 
we can explicitly construct all pfaffian representations of $C$ from a given one (for example, from a decomposable
representation with symmetric blocks induced by one of the $2^{g-1}(2^{g}+1)$ even theta characteristics). 
Section~\ref{seclast} is an exposition on plane quartics.
Concrete examples and algorithms for computing Aronhold bundles and theta characteristics are given.
 
\section{The Canonical Form} 
\label{canform} 

Canonical forms play an important role in explicit descriptions of the moduli space $M_C(2,K_C)$. 
In the sequel we outline an algorithm for such computation. 

Let $A=x_0 A_0+x_2 A_2+x_2 A_2$ be a representation of $C$. 
We can always assume that after a projective change of coordinates $C$ intersects the line $L: x_0=0$ in distinct points
$P_1=(p_1,1,0),\ldots,P_d=(p_d,1,0)$. 
We will prove that $A$ is equivalent to a representation in the canonical form.

\begin{proposition}\label{propcanform}
For every pfaffian representation $A$ of $C$ there exists a basis of $k^{2d}$ in which $A$  
has the canonical form
 \begin{equation}
\label{caneqfor}A=x_1\left[\begin{array}{cccc}
 I & 0   &   \cdots & 0   \\
 0 & I   &   \cdots & 0   \\
 \vdots  & &   \ddots & \vdots  \\
 0 &  0    &  \cdots  & I
 \end{array}\right]-x_2\left[\begin{array}{cccc}
 D_1 & 0   &   \cdots & 0   \\
 0 & D_2   &   \cdots & 0   \\
 \vdots  & &   \ddots & \vdots  \\
 0 &  0    &  \cdots  & D_d
 \end{array}\right]+x_0 A_0,\end{equation}
 where 
 $$I=\left[\begin{array}{cc}
 0 &  1   \\
 -1  & 0  
 \end{array}\right]\ \mbox{ and }\ D_i=\left[\begin{array}{cc}
 0 &  p_i   \\
 -p_i  & 0  
\end{array}\right].$$
\end{proposition}
\begin{proof} As above, let
$P_i=(p_i,1,0)$ be $d$ distinct points in the intersection of $F(x_0,x_1,x_2)=0$ with the line $L: x_0=0$.
By restricting to $L$, we obtain the pencil of skew-symmetric matrices
$$A_L=x_1A_1+x_2A_2$$
with $\pf A_L=F|_L=F(0,x_1,x_2)=\prod_{i=1}^d(x_1-p_i x_2)$.
Since all the matrices in the representation are skew-symmetric, $\coker A$ and $\ker A$ define the same rank 2 vector bundle $\cE$.
Note that $\cE(P_i)$ is the kernel of $p_i A_1+A_2$.
Thus $\cE(P_i),\ i=1,\ldots,d$ are $2$-dimensional subspaces in $k^{2d}$. 
In the proof of~\cite[Proposition 3.4]{anitapfaff} we showed that $h^0(C,\cE(-1))=0$ implies the following important fact:
the union of bases of the vector spaces $\cE(P_i)$
span the whole space $k^{2d}$.
In this basis $A_L$  is equivalent to
$$J=\left[\begin{array}{cccc}
J_1 & 0   &   \cdots & 0   \\
0 & J_2   &   \cdots & 0   \\
\vdots  & &   \ddots & \vdots  \\
0 &  0    &  \cdots  & J_d
\end{array}\right]
$$
where 
$$J_i=\left[\begin{array}{cc}
0 &  x_1-p_i x_2   \\
-x_1+p_i x_2 & 0  
\end{array}\right].$$ No $2\times 2$ block $J_i$ is identically equal to $0$,
otherwise the rank of $A$ at $P_i$ would be at most $2d-4$.
\end{proof}

Let $A$ be a representation of $C$ and $\cE$ its cokernel. Choose the basis of $k^{2d}$
$$\cup_{i=1}^d \left\{\mbox{2--dimensional basis of }\cE(P_i)\right\}.$$ 
and denote by $P\in\GL(2d,k)$ the change of basis matrix. Then $P$ acts on the representation by $A\mapsto P\, A\, P^t$.
Thus $A$ is equivalent to a representation in the canonical form by Proposition~\ref{propcanform}.

In~\cite{anitapfaff} we established a one to one correspondence between the pfaffian representations (up to equivalence) of $C$ 
and the open set  
\begin{equation} \label{moduliopsub} M_C(2,K_C)-\{\mathcal{K}\, :\, h^0(C, \mathcal{K})>0 \}.\end{equation}
It thus suffices to find all pfaffian representations of $C$ in 
the canonical form, which will yield a set of equations describing the open subset~(\ref{moduliopsub}) in the moduli space.

There are $d(2d-1)$ parameters in the representation (\ref{caneqfor}), namely the entries of the $2d\times 2d$ skew-symmetric matrix
$A_0=[a^0_{ij}]$.
Since $\pf A$ equals $F$, we get $\frac{1}{2}d(d+1)$ relations among $a^0_{ij}$. Indeed, every monomial 
$x_0^{\alpha} x_1^{\beta} x_2^{\gamma},\ \alpha\neq 0,\, \alpha+\beta+\gamma=d$ in $F$ gives one equation.
By the implicit function theorem we are left with $d(2d-1)-\frac{1}{2}d(d+1)=\frac{3}{2}d(d-1)$ parameters $a^0_{ij}$.
Recall that pfaffian representations are equivalent under the action
$$A\ \mapsto\ R\, A\, R^t,\ \ R\in\GL(2d,k).$$
By a suitable $R$ we can further reduce the number of parameters in $A_0$. 
Of course we only consider $R$ whose action preserves the canonical form.
This way we will reduce the number of equivalent representations in each equivalence class to 1.
The following lemma is an elementary exercise in linear algebra.
\begin{lemma} \label{lemreducepf}
The action
$A\ \mapsto\ R\, A\, R^t$ preserves the canonical form of the first two matrices in the representation if and only if  
$R$ equals 
$$\left[\begin{array}{cccc}
R_1 & 0   &   \cdots & 0   \\
0 & R_2   &          & 0   \\
\vdots  & &   \ddots & \vdots  \\
0 &  0    &  \cdots  & R_d
\end{array}\right],
$$ 
where every $R_i$ is an invertible $2\times 2$ matrix with determinant 1.
\end{lemma}
Matrix $R$ in Lemma~\ref{lemreducepf} acts on the third matrix in the representation by
\begin{equation*}R\, A_0\, R^t= \left[\begin{array}{ccl}
0 &       &         \\
  &\ddots &      R_i\, \left[\begin{array}{cc}
a^0_{2i-1,2j-1} &  a^0_{2i-1,2j}   \\
a^0_{2i,2j-1} & a^0_{2i,2j} 
\end{array}\right] \, R_j^t   \\
  &       &   \ddots                \\
  &       &  \ \ \ \ \ \ \ \ 0         
\end{array}\right],\ \ i,j=1,\ldots,d.\end{equation*}
 In other words, view $A_0$ as $d\times d$ matrix of $2\times 2$ blocks. Then $R$ acts on $i,j-$th block by 
$\ast \mapsto R_i\, \ast \, R_j^t$. 
In every $R_i$ we have a choice of 3 independent parameters since its determinant is 1. 
Thus each $R_i$ reduces the number of $a^0_{ij}$'s by 3.

For a general $F(x_0,x_1,x_2)$ this sequence of reductions can be performed explicitly. We are left with 
\begin{enumerate}
\item[] $3(g-1)=\frac{3}{2}d(d-3)=$
\item[] $d(2d-1)$ [number of parameters $a^0_{ij}$ in $A$ in the canonical form] $-$ 
\item[] $\frac{1}{2}d(d+1)$ [relations since $\pf A=F$]$-$
\item[] $3 d$ [parameters reduced by the action of equivalence relation $R\cdot A\cdot R^t$]
\end{enumerate}
independent variables $a^0_{ij}$. As expected, this represents the open subset~(\ref{moduliopsub}) of the moduli space in $\PP^{d(2d-1)-1}$.
\begin{remark} \label{decremex} \rm{
In the above considerations we could take any other canonical form of the matrix pair $A_1,A_2$. For example, 
the equivalence relation action $Q\, A\, Q^t$ of
 $$Q=\left[\begin{array}{cccccc}
1 & -1 & 0 & 0 &  \cdots & 0   \\
0 &  0 & 1 &-1 &       & 0   \\
\vdots& &  & \ddots & & \vdots  \\
0 &  0  &    \cdots & & 1 & -1\\
0 & 1 & 0 & 0 &  \cdots & 0   \\
0 &  0 & 0 &1 &       & 0   \\
\vdots& &  & \ddots & & \vdots  \\
0 &  0  &    \cdots &  & 0 & 1
\end{array}\right]
$$ 
 brings the first two matrices in~(\ref{caneqfor}) into
\begin{equation}\label{seccanform}
\left[ \begin{array}{cc}
0 & \Id \\
-\Id & 0
\end{array}\right], \ \ \left[ \begin{array}{cc}
0 & D \\
-D & 0
\end{array}\right],
\end{equation}
where $D$ is the diagonal matrix $\{p_1,\ldots,p_d\}$.
 This canonical form is particularly useful since it naturally includes all the decomposable representations~(\ref{decomrep}). 
The same canonical form was obtained in~\cite{rodman} , where canonical forms for matrix pairs were classified purely by the methods of
linear algebra. 
}\end{remark}


\section{Tangents} 
\label{tangents}

In this section we explain how tangents or lines through two points $\lambda,\mu$ on $C$ can be read from pfaffian representations. Moreover,
for any pfaffian representation $A$ we relate vectors in $\coker A(\lambda)$ and $\coker A(\mu)$.

Since our aim is to do concrete calculations, we describe the correspondence in Theorem~\ref{beauvilcor} explicitly:
Denote by $\pf^{ij} A$ the pfaffian of the $(2d-2)\times(2d-2)$ skew-symmetric matrix obtained 
by removing the $i$th and $j$th rows and columns from $A$. It is easy to prove the following analogue of the
Jacobi's formula for the derivative of determinants:
\begin{equation}\label{odvod}\frac{\partial F(x_0,x_1,x_2)}{\partial x_k}= \sum_{i,j} a_{ij}^k\ (-1)^{i+j}\pf^{ij} A(x_0,x_1,x_2).\end{equation}
If for some $x=(x_0,x_1,x_2)\in C$ all $2d-2$ pfaffian minors vanish, then $x$ must be a singular point of $F$.
By our assumption $F$ is smooth, thus $\rank A(x) \geq 2d-2$ for all $x\in C$.
Rank of skew-symmetric matrices is even and $\det A=F^2=0,$ thus equality holds. Therefore $\coker A(x)$ defines
a $\rank 2$ bundle over $C$.
Define the \textit{pfaffian adjoint} of $A$ to be the skew-symmetric matrix
\begin{equation*}\tilde{A}= \left[\begin{array}{ccl}
0 &       &         \\
  &\ddots &      (-1)^{i+j} \pf^{ij} A   \\
  &       &   \ddots                \\
  &       &  \ \ \ \ \ \ \ \ 0         
\end{array}\right].\end{equation*}
Again, by analogy with determinants, the following holds
\begin{equation} \label{pfaflinalg} 
\tilde{A} \, A=\pf \! A\ \mbox{Id}_{2d}.\end{equation}
More properties and linear algebra of pfaffians can be found in~\cite[Appendix D]{fulton}.

Since $\pf \! A=c\, F$, the cokernel 
can be obtained from $\tilde{A}$ by using (\ref{pfaflinalg}). Indeed, for any point $x\in C$ 
every column (or row) of $\tilde{A}(x)$ is in $\coker A(x)$.
Since $A(x)$ is a skew-symmetric matrix with linear entries, we have
$$\coker A(x)\cong\ker A(x)$$ which we denote  by $\cE(x)$. We will often view $\cE$ as a $\PP^1$--bundle,
or equivalently as a ruled surface $\PP\cE$.

\begin{lemma} \label{lemtang}
Every representation of $C$ yields tangents and lines through $\lambda,\mu\in C$ in the following way:
\begin{itemize}
\item Let $v_{\lambda}\in \cE(\lambda),\, u_{\mu}\in \cE(\mu)$ be independent vectors. Then
$u_{\mu}^t A(x) v_{\lambda}$ is either identically $0$ or defines a line through ${\lambda}$ and $\mu$.
\item Let $\mathcal{L}in\{u_{\lambda}, v_{\lambda}\}=\cE(\lambda)$. Then 
$u_{\lambda}^tA(x)v_{\lambda}=0$ is an equation of the tangent line at ${\lambda}\in C$.
\item Analogously, for $\mathcal{L}in\{u_{\mu}, v_{\mu}\}= \cE(\mu)$ the equation 
$u_{\mu}^tA(x)v_{\mu}=0$ defines the tangent at ${\mu}\in C$.
\end{itemize}
\end{lemma}
\begin{proof} The first assertion is obvious as $u_{\mu}^tA(x)v_{\lambda}$ is linear and  
equals $0$ at the points $\lambda=(\lambda_0,\lambda_1,\lambda_2)$ and $\mu=(\mu_0,\mu_1,\mu_2)$.

By~(\ref{pfaflinalg}), $\tilde{A}$ is a rank 2 skew-symmetric matrix at the points of $C$.
Using the Pl\"{u}cker coordinates in~(\ref{explrk}),  $u_{\lambda}\wedge v_{\lambda}$  equals to a multiple of $\tilde{A}(\lambda)$.
Then by~(\ref{odvod})
\begin{eqnarray*}
\sum_{k=0}^2 x_k\, \frac{\partial F}{\partial x_k}(\lambda)&=&\sum_{k=0}^2 x_k \sum_{i,j} a_{ij}^k \ (-1)^{i+j}\pf^{ij} A(\lambda)\\
  &=& \sum_{i,j} (a_{ij}^0 x_0+a_{ij}^1 x_1+a_{ij}^2 x_2)\ (-1)^{i+j}\pf^{ij} A(\lambda)\\
  &=& \trace\left[A(x)\cdot \tilde{A}(\lambda)\right],
\end{eqnarray*}
which equals to a nonzero multiple of 
$\trace\left[A(x)\cdot (u_{\lambda}\wedge v_{\lambda})\right]=2\ u_{\lambda}^tA(x)v_{\lambda}.$
In particular we proved that $u_{\lambda}^tA(x)v_{\lambda}$ is not identically $0$.

The same way we obtain the tangent to $C$ at $\mu$.
\end{proof}

$\cE(x)$ can be viewed as pairs of independent vectors $[u_x,v_x]$ modulo right--hand--side multiplications by invertible $2\times 2$ 
matrices. Then $T,T'\in \GL(2,k)$ can be chosen so that 
$$T\, [u_{\lambda}, v_{\lambda}]^t\, A(x)\, [u_{\mu}, v_{\mu}] \, T'=
\left[ \begin{array}{cc} u_{\lambda}^t A(x) u_{\mu} & u_{\lambda}^t A(x) v_{\mu} \\ 
v_{\lambda}^t A(x) u_{\mu}& v_{\lambda}^t A(x) v_{\mu} \end{array} \right]$$
equals either 
$$F(x) \left[ \begin{array}{cc} 0&0\\ 0& 0 \end{array} \right],\ F(x) \left[ \begin{array}{cc} 0&0\\ 0& 1 \end{array} \right]
\mbox{ or } F(x) \left[ \begin{array}{cc} 0&1\\ 1& 0 \end{array} \right].$$
This enables us to relate points on $C$.
A distinct pair of points $\lambda, \mu$ on $C$ with respect to representation $A(x)$ is either: 
\begin{itemize}
\item \textit{inadmissible} if $u_{\mu}^t A(x) v_{\lambda}\equiv0$ on $\PP^2$ for all $v_{\lambda}\in \PP\cE(\lambda),\, u_{\mu}\in \PP\cE(\mu)$;
\item \textit{semiadmissible} if there exists a unique pair of vectors $u_{\lambda}\in \PP\cE(\lambda),\, u_{\mu}\in \PP\cE(\mu)$ such that
$w_{\mu}^t A(x) u_{\lambda}\equiv 0$ and $u_{\mu}^t A(x) w_{\lambda}\equiv 0$  for all $w_{\lambda}\in \PP\cE(\lambda),\, w_{\mu}\in \PP\cE(\mu)$;
\item \textit{admissible} if for every $v_{\lambda}\in \PP\cE(\lambda)$ there 
exists exactly one $ v_{\mu}\in \PP\cE(\mu)$
for which $v_{\mu}^t A(x) v_{\lambda}\equiv 0$ on the whole $\PP^2$. 
\end{itemize}
For admissible $\lambda, \mu$ we thus obtain a one to one correspondence between vectors in $\PP\cE(\lambda)$ and $\PP\cE(\mu)$. Namely,
$ v_{\lambda}$ corresponds to $v_{\mu}$ iff $v_{\mu}^t A(x) v_{\lambda}\equiv 0$. Consider next a semiadmissible pair $\lambda, \mu\in C$  and 
$v_{\lambda}\in \PP\cE(\lambda),\, v_{\mu}\in \PP\cE(\mu)$. We claim that $v_{\mu}^t A(x) v_{\lambda}\equiv 0$ if and only if either 
$v_{\lambda}=u_{\lambda}$ or $v_{\mu}=u_{\mu}$. Indeed, assume that $v_{\lambda}\neq u_{\lambda}$, then by linearity 
$v_{\mu}^t A(x) w_{\lambda}\equiv 0$ for any 
$w_{\lambda}\in \mathcal{L}in\{ v_{\lambda},u_{\lambda} \}=\cE(\lambda)$. Thus by uniqueness in the definition $v_{\mu}=u_{\mu}$. 

In the proof of Lemma~\ref{lemtang} we showed 
\begin{corollary} A bitangent line through $\lambda, \mu\in C$ is defined by either $u_{\mu}^t A(x) v_{\lambda}$, $u_{\lambda}^t A(x) v_{\lambda}$
or $u_{\mu}^t A(x) v_{\mu}$.
\end{corollary}
\begin{corollary}
Let  $u_{\lambda},v_{\lambda}\in \PP\cE(\lambda)$. Then $u_{\lambda}^t A(x) v_{\lambda}\equiv 0$ if and only if $u_{\lambda}=v_{\lambda}$. 
\end{corollary}
This means that $\lambda,\lambda$ can be always considered  as an admissible pair of points on $C$.

\section{Elementary Transformations} 
\label{beginsec}

In this section we define \textit{elementary transformations of pfaffian representations} of the curve $C:\, F(x_0,x_1,x_2)=0$ in $\PP^2$, 
which is not necessarily smooth or irreducible. The standard notation of vessels~\cite{vinn2} will be used.

Choose the coordinates of $F$ so that the line $\{x_0=0\}$ and the curve intersect in $d$ distinct smooth points.
For the sake of clearer notation we move to the affine plane.
Consider a linear pfaffian representation
$$ \pf (y_1 \sigma_2-y_2 \sigma_1+\gamma)=c\, f(y_1,y_2)\ \ \ \ \ \ c\neq 0,$$
where $y_1=\frac{x_1}{x_0},\, y_1=\frac{x_1}{x_0}$ are the affine coordinates $f(y_1,y_2)=F(1,y_1,y_2)$ and
$\sigma_1,\sigma_2, \gamma$ are $2d\times 2d$ skew-symmetric matrices.
Denote $$\coker (y_1 \sigma_2-y_2 \sigma_1+\gamma)\cong\ker (y_1 \sigma_2-y_2 \sigma_1+\gamma)$$ by $\cE(y_1,y_2)$. 

Pick distinct regular affine points $\lambda=(\lambda_1,\lambda_2)$ and $\mu=(\mu_1,\mu_2)$ on $C$. Then each
$\cE(\lambda),\, \cE(\mu)$ is a 2 dimensional vector space in $k^{2d}$. 
For all $v_{\lambda}\in \cE(\lambda),\, u_{\mu}\in \cE(\mu)$ 
\begin{equation}
\label{vin9}
\begin{array}{c}
(\lambda_1 \sigma_2-\lambda_2 \sigma_1+\gamma)v_{\lambda}=0 \\
(\mu_1 \sigma_2-\mu_2 \sigma_1+\gamma)u_{\mu}=0
\end{array}
\end{equation}
and 
\begin{equation}
\label{vin9a}
\begin{array}{c}
(y_1 \sigma_2-y_2 \sigma_1+\gamma)v_{\lambda}=((y_1-\lambda_1) \sigma_2-(y_2-\lambda_2) \sigma_1)v_{\lambda} \\
(y_1 \sigma_2-y_2 \sigma_1+\gamma)u_{\mu}=\left((y_1-\mu_1) \sigma_2-(y_2-\mu_2) \sigma_1\right)u_{\mu}
\end{array}
\end{equation}
hold. 
Thus 
\begin{eqnarray} \label{vin10}
v_{\lambda}^t(\lambda_1 \sigma_2-\lambda_2 \sigma_1+\gamma)u_{\mu}=0\ \mbox{ and }\ v_{\lambda}^t(\mu_1 \sigma_2-\mu_2 \sigma_1+\gamma)u_{\mu}=0,
\end{eqnarray}
which implies
\begin{eqnarray*}
(\lambda_1-\mu_1)v_{\lambda}^t \sigma_2 u_{\mu}= (\lambda_2-\mu_2)v_{\lambda}^t \sigma_1 u_{\mu}.
\end{eqnarray*}
In other words, for any pair of complex parameters $t_1,t_2$, 
\begin{eqnarray}
\label{vin6.7}
\frac{1}{t_1(\lambda_1-\mu_1)+t_2(\lambda_2-\mu_2)}\, v_{\lambda}^t (t_1 \sigma_1 +t_2 \sigma_2) u_{\mu}
\end{eqnarray}
is constant whenever the denominator is nonzero. Denote this constant by $K_{v_{\lambda}u_{\mu}}$.

The pair of vectors $v_{\lambda}\in \cE(\lambda),\, u_{\mu}\in \cE(\mu)$ is called \textit{admissible} if $K_{v_{\lambda}u_{\mu}}$ is not 0.
Define 
\begin{eqnarray*}
\label{vin12}
\tilde{\gamma} &=&\gamma+\frac{1}{K_{v_{\lambda}u_{\mu}}}\left(\sigma_1(u_{\mu}v_{\lambda}^t+v_{\lambda}u_{\mu}^t)\sigma_2-\sigma_2(u_{\mu}v_{\lambda}t+v_{\lambda}u_{\mu}^t)\sigma_1 \right)\\
               &=&\gamma+\frac{1}{K_{v_{\lambda}u_{\mu}}}\left(\sigma_1(u_{\mu}v_{\lambda}^t+v_{\lambda}u_{\mu}^t)\sigma_2-(\sigma_1(u_{\mu}v_{\lambda}^t+v_{\lambda}u_{\mu}^t)\sigma_2)^t\right),
\end{eqnarray*}
and
\begin{eqnarray*}
\label{vin12a}
\bar{\gamma} &=&\gamma+2 \rho \left(\sigma_1v_{\lambda} v_{\lambda}^t\sigma_2-\sigma_2 v_{\lambda}v_{\lambda}t\sigma_1 \right)\\
               &=&\gamma+ 2\rho \left( \sigma_1 v_{\lambda}v_{\lambda}^t\sigma_2-(\sigma_1v_{\lambda}v_{\lambda}^t\sigma_2)^t \right),
\end{eqnarray*}
where $u_{\mu}v_{\lambda}^t,\ v_{\lambda}u_{\mu}^t$ and $v_{\lambda}v_{\lambda}^t$ are rank 1 matrices and $\rho$ is an arbitrary constant.
It is obvious that $\tilde{\gamma}=-\tilde{\gamma}^t$ and $\bar{\gamma}=-\bar{\gamma}^t$ are skew-symmetric.\\

In the sequel the following properties of the wedge product will be useful: for any skew-symmetric matrix  $\sigma$
and vectors $w,v$ holds
$(\sigma w)^t=-w^t \sigma,\ w^t \sigma v=-v^t \sigma w$ and $w\wedge v=2(wv^t-vw^t)$.

Using the above, $\tilde{\gamma}$ rewrites into 
\begin{eqnarray}
\label{tilgam}
\tilde{\gamma} =\gamma-\frac{1}{2\, K_{v_{\lambda}u_{\mu}}}\, \sigma_1u_{\mu} \wedge \sigma_2 v_{\lambda}+
\frac{1}{2\, K_{v_{\lambda}u_{\mu}}}\, \sigma_2 u_{\mu} \wedge \sigma_1 v_{\lambda}
\end{eqnarray}
and $\bar{\gamma}$ rewrites into 
\begin{eqnarray}
\label{bargam}
\bar{\gamma} =\gamma+
 \rho\ \sigma_2 v_{\lambda} \wedge \sigma_1 v_{\lambda}.
\end{eqnarray}

\begin{theorem}
Let $y_1 \sigma_2-y_2 \sigma_1+\gamma$ be a representation of $C$ and $\cE(y_1,y_2)$ its cokernel. Choose 
$v_{\lambda}\in \cE(\lambda),\, u_{\mu}\in \cE(\mu)$ such that $K_{v_{\lambda}u_{\mu}}\neq 0$. Then
$y_1 \sigma_2-y_2 \sigma_1+\tilde{\gamma}$ and $y_1 \sigma_2-y_2 \sigma_1+\bar{\gamma}$ are pfaffian representations of $C$
since
$$\pf (y_1 \sigma_2-y_2 \sigma_1+\tilde{\gamma})=\pf (y_1 \sigma_2-y_2 \sigma_1+\gamma)=\pf (y_1 \sigma_2-y_2 \sigma_1+\bar{\gamma}).$$
\end{theorem}

\begin{proof}
For any $T$ and skew-symmetric $S$ holds 
$$\pf (T S T^t)=\det T \, \pf S$$ 
by~\cite[Appendix D]{fulton}. Recall also that for any two square matrices $T',T''$
$$\det(\Id+T' T'')=\det(\Id+T'' T').$$
Therefore proving the equality
\begin{equation} \label{vin14} \begin{array}{c}
\left(\Id+ \frac{(t_1 \sigma_1 +t_2 \sigma_2)u_{\mu}v_{\lambda}^t}{K_{v_{\lambda}u_{\mu}}\left(t_1(y_1-\lambda_1)+t_2(y_2-\lambda_2)\right)} \right)(y_1 \sigma_2-y_2 \sigma_1+\gamma)
\left(\Id-\frac{v_{\lambda}u_{\mu}^t(t_1 \sigma_1 +t_2 \sigma_2)}{K_{v_{\lambda}u_{\mu}}\left(t_1(y_1-\lambda_1)+t_2(y_2-\lambda_2)\right)}  \right)\\
 \parallel  \\
\left(\Id- \frac{(t_1 \sigma_1 +t_2 \sigma_2) v_{\lambda}u_{\mu}^t}{K_{v_{\lambda}u_{\mu}}\left(t_1(y_1-\lambda_1)+t_2(y_2-\lambda_2)\right)}\right)(y_1 \sigma_2-y_2 \sigma_1+\tilde{\gamma})
\left(\Id+ \frac{u_{\mu}v_{\lambda}^t(t_1 \sigma_1 +t_2 \sigma_2)}{K_{v_{\lambda}u_{\mu}}\left(t_1(y_1-\lambda_1)+t_2(y_2-\lambda_2)\right)} \right) 
\end{array}\end{equation}
will imply the first statement. (\ref{vin14}) equals
\begin{eqnarray*}
y_1 \sigma_2-y_2 \sigma_1+\gamma+ \frac{1}{2K_{v_{\lambda}u_{\mu}}} (y_1 \sigma_2-y_2 \sigma_1+\gamma)v_{\lambda} \wedge 
\frac{t_1 \sigma_1 +t_2 \sigma_2}{t_1(y_1-\lambda_1)+t_2(y_2-\lambda_2)} u_{\mu}\ &= \\
y_1 \sigma_2-y_2 \sigma_1+\tilde{\gamma}- \frac{1}{2K_{v_{\lambda}u_{\mu}}}(y_1 \sigma_2-y_2 \sigma_1+\tilde{\gamma})u_{\mu} \wedge 
\frac{t_1 \sigma_1 +t_2 \sigma_2}{t_1(y_1-\lambda_1)+t_2(y_2-\lambda_2)} v_{\lambda}. &
\end{eqnarray*}

From~(\ref{vin6.7}) we obtain that 
\begin{eqnarray*}
-\frac{1}{2K_{v_{\lambda}u_{\mu}}}\, \left(\sigma_1u_{\mu} \wedge \sigma_2 v_{\lambda}-
 \sigma_2 u_{\mu} \wedge \sigma_1 v_{\lambda}\right) u_{\mu} &=\\
 \frac{1}{K_{v_{\lambda}u_{\mu}}} \left( \sigma_1u_{\mu} v^t_{\lambda}\sigma_2- \sigma_2 v_{\lambda}u^t_{\mu} \sigma_1- 
 \sigma_2u_{\mu} v^t_{\lambda}\sigma_1 +\sigma_1 v_{\lambda}u^t_{\mu} \sigma_2 \right)u_{\mu} & =\\
  \frac{1}{K_{v_{\lambda}u_{\mu}}} \left( \sigma_1u_{\mu}\, v^t_{\lambda}\sigma_2u_{\mu}- 
 \sigma_2u_{\mu}\, v^t_{\lambda}\sigma_1u_{\mu}  \right)=(\lambda_2-\mu_2) \sigma_1u_{\mu}-(\lambda_1-\mu_1) \sigma_2u_{\mu}.&
\end{eqnarray*} 
Analogously 
\begin{eqnarray*}
\frac{1}{2K_{v_{\lambda}u_{\mu}}}\, \left(\sigma_1u_{\mu} \wedge \sigma_2 v_{\lambda}-
 \sigma_2 u_{\mu} \wedge \sigma_1 v_{\lambda}\right)v_{\lambda}=(\lambda_2-\mu_2) \sigma_1 v_{\lambda}-(\lambda_1-\mu_1) \sigma_2 v_{\lambda}.
\end{eqnarray*} 
Together with~(\ref{vin9a}) this implies
\begin{equation}
\label{vin13}
\begin{array}{c}  
(\lambda_1 \sigma_2-\lambda_2 \sigma_1+\tilde{\gamma})u_{\mu}=0 \\
(\mu_1 \sigma_2-\mu_2 \sigma_1+\tilde{\gamma})v_{\lambda}=0
\end{array}
\end{equation}
and
\begin{equation}
\label{vin13a}
\begin{array}{c}
(y_1 \sigma_2-y_2 \sigma_1+\tilde{\gamma})v_{\lambda}=((y_1-\mu_1) \sigma_2-(y_2-\mu_2) \sigma_1)v_{\lambda} \\
(y_1 \sigma_2-y_2 \sigma_1+\tilde{\gamma})u_{\mu}=\left((y_1-\lambda_1) \sigma_2-(y_2-\lambda_2) \sigma_1\right)u_{\mu}
\end{array}
\end{equation}
In other words, $u_{\mu}\in \tilde{\cE}(\lambda)$ and $v_{\lambda}\in \tilde{\cE}(\mu)$, where 
$\tilde{\cE}(y_1,y_2)=\coker (y_1 \sigma_2-y_2 \sigma_1+\tilde{\gamma})$. 

Now it is easy to verify that 
\begin{eqnarray*}
((y_1-\lambda_1) \sigma_2-(y_2-\lambda_2) \sigma_1)v_{\lambda} \wedge 
\frac{t_1 \sigma_1 +t_2 \sigma_2}{t_1(y_1-\lambda_1)+t_2(y_2-\lambda_2)} u_{\mu}\ &= \\
 -\sigma_1u_{\mu} \wedge \sigma_2 v_{\lambda}+\sigma_2 u_{\mu} \wedge \sigma_1 v_{\lambda}-
 \ \ \ \ \ \ \ \ \ \ \ \ \ \ \ \ \ \ \ \ \ \ \ \ \ \ \ \ \ \ \ \ \ \ \ \ \ \ \  &\\
 ((y_1-\lambda_1) \sigma_2-(y_2-\lambda_2) \sigma_1)u_{\mu} \wedge 
\frac{t_1 \sigma_1 +t_2 \sigma_2}{t_1(y_1-\lambda_1)+t_2(y_2-\lambda_2)} v_{\lambda}.& 
\end{eqnarray*}
which together with~(\ref{vin9a}) and~(\ref{vin13a}) finishes the proof of~(\ref{vin14}).\\

The second statement will be proved if we show
\begin{equation}
\label{vin14sec}
\begin{array}{c}
 (y_1 \sigma_2-y_2 \sigma_1+\gamma)\\
\parallel  \\
\left(\Id-2 \rho \frac{t_1 \sigma_1 +t_2 \sigma_2}{t_1(y_1-\lambda_1)+t_2(y_2-\lambda_2)} v_{\lambda}v_{\lambda}^t \right)
(y_1 \sigma_2-y_2 \sigma_1+\bar{\gamma})
\left(\Id+2 \rho v_{\lambda} v_{\lambda}^t \frac{t_1 \sigma_1 +t_2 \sigma_2}{t_1(y_1-\lambda_1)+t_2(y_2-\lambda_2)} \right) \\
\parallel  \\
y_1 \sigma_2-y_2 \sigma_1+\bar{\gamma}-\rho\ (y_1 \sigma_2-y_2 \sigma_1+\tilde{\gamma})v_{\lambda} \wedge 
\frac{t_1 \sigma_1 +t_2 \sigma_2}{t_1(y_1-\lambda_1)+t_2(y_2-\lambda_2)} v_{\lambda}. 
\end{array}
\end{equation}
Note that
\begin{eqnarray*}
\left( \sigma_2 v_{\lambda} \wedge \sigma_1 v_{\lambda}\right) v_{\lambda} &=&\\
  \left( - \sigma_2 v_{\lambda}v^t_{\lambda} \sigma_1+\sigma_1 v_{\lambda}v^t_{\lambda} \sigma_2 \right)v_{\lambda} & =&0
\end{eqnarray*} 
implies
\begin{eqnarray}\label{vintype2}
(\lambda_1 \sigma_2-\lambda_2 \sigma_1+\bar{\gamma})v_{\lambda}&=&0\\
(y_1 \sigma_2-y_2 \sigma_1+\bar{\gamma})v_{\lambda}&=&((y_1-\lambda_1) \sigma_2-(y_2-\lambda_2) \sigma_1)v_{\lambda} .\nonumber
\end{eqnarray}
This means that $v_{\lambda}\in \bar{\cE}(\lambda)$, where  $\bar{\cE}(y_1,y_2)$ denotes  the cokernel of 
$y_1 \sigma_2-y_2 \sigma_1+\bar{\gamma}$. 

The above and the obvious equality 
\begin{eqnarray*}
 -\sigma_2 v_{\lambda} \wedge \sigma_1 v_{\lambda}=
 \ \ \ \ \ \ \ \ \ \ \ \ \ \ \ \ \ \ \ \ \ \ \ \ \ \ \ \ \ \ \ \ \ \ \ \ \ \ \ \ \ \ \  &\\
 -((y_1-\lambda_1) \sigma_2-(y_2-\lambda_2) \sigma_1)v_{\lambda} \wedge 
\frac{t_1 \sigma_1 +t_2 \sigma_2}{t_1(y_1-\lambda_1)+t_2(y_2-\lambda_2)} v_{\lambda}.& 
\end{eqnarray*}
finish the proof of~(\ref{vin14sec}).
\end{proof}

By performing \textit{elementary transformations} of $y_1 \sigma_2-y_2 \sigma_1+\gamma$,
we constructed new pfaffian representations of $C$:
\begin{definition}\label{defelemrepr}\rm{
The \textit{Type I elementary transformation}  $y_1 \sigma_2-y_2 \sigma_1+\tilde{\gamma}$ based on the admissible vectors 
$v_{\lambda}\in \cE(\lambda),\, u_{\mu}\in \cE(\mu)$, \\
The \textit{Type II elementary transformation} $y_1 \sigma_2-y_2 \sigma_1+\bar{\gamma}$ based on $v_{\lambda}\in \cE(\lambda)$ and the 
constant $\rho\neq 0$.}
\end{definition}
 
The fact that $v_{\lambda}\in \tilde{\cE}(\mu),\ u_{\mu}\in \tilde{\cE}(\lambda)$ and 
$v_{\lambda}\in \bar{\cE}(\lambda)$ implies the following
\begin{corollary} The Type I elementary transformation of 
$y_1 \sigma_2-y_2 \sigma_1+\tilde{\gamma}$ 
based on $u_{\mu}\in \tilde{\cE}(\lambda),\ v_{\lambda}\in \tilde{\cE}(\mu)$ brings us back to $y_1 \sigma_2-y_2 \sigma_1+\gamma$.
The same way the Type II elementary transformation of 
$y_1 \sigma_2-y_2 \sigma_1+\bar{\gamma}$ 
based on $v_{\lambda}\in \bar{\cE}(\lambda)$ and $-\rho$ brings us back to $y_1 \sigma_2-y_2 \sigma_1+\gamma$.
\end{corollary}


It can be easily seen that the Type I and II elementary transformations are special rank 2 cases of 
"the concrete interpolation problem for meromorphic bundle maps" studied in~\cite{vinn2}:
For a given array of $m$ constants $\rho_i\in k$, $m$ distinct points $\lambda^i\in C$ and  vectors $w_{i}\in \cE(\lambda^i)$ 
define the \textit{Type CONINT} elementary transformation of $y_1 \sigma_2-y_2 \sigma_1+{\gamma}$ by
$$\check{\gamma}=\gamma+\sigma_1 w \Gamma^{-1}w^t \sigma_2-\sigma_2 w \Gamma^{-1}w^t \sigma_1,$$
where $w=[w_1,\ldots,w_m]$ and
$-\Gamma$ equals the $m\times m$ symmetric matrix with $\rho_i$'s along the diagonal and $K_{w_{i}\, w_{j}}$ 
at the $(i<j)$--th position. Then $y_1 \sigma_2-y_2 \sigma_1+\check{\gamma}$ is indeed a pfaffian representation of $C$
since 
$y_1 \sigma_2-y_2 \sigma_1+{\gamma}=Z^t(y)\, (y_1 \sigma_2-y_2 \sigma_1+\check{\gamma})\, Z(y)$ for 
$$Z(y)=\Id+w\, \mbox{Diagonal}\! \left[t_1(y_1-\lambda_1^i)+t_2(y_2-\lambda_2^i)\right]^{-1}_{\substack{i=1\ldots m}}\, 
\Gamma^{-1}\,w^t\, (t_1 \sigma_1+t_2 \sigma_2) .$$
An easy exercise in linear algebra shows that 
$$Z^{-1}(y)=\Id-w\, \Gamma^{-1}\, \mbox{Diagonal}\! \left[t_1(y_1-\lambda_1^i)+t_2(y_2-\lambda_2^i)\right]^{-1}_{\substack{i=1\ldots m}}\, 
w^t\, (t_1 \sigma_1+t_2 \sigma_2) $$
and $\det Z(y)=\det Z^{-1}(y)=1$.
Observe that the only condition for $\check{\gamma}$ to exist is that $\Gamma$ is invertible.

In Section~\ref{decsec} we will show that any
two pfaffian representations of $C$ can be bridged by a finite sequence of Type I and Type II elementary transformations.
Thus it is enough to restrict our study to these two types.

\section{Admissible arrays}
\label{sectfour}
 
We briefly relate the definitions of admissible pairs of points and admissible pairs of vectors. Let 
\begin{eqnarray*}A(x_0,x_1,x_2)&=&x_1 A_1+x_2 A_2+x_0 A_0\\
&=&x_1 \sigma_2-x_2 \sigma_1+x_0 \gamma
\end{eqnarray*}
be a representation of $C$ with 
the cokernel $\cE(x)$. As before, let
$v_{\lambda}\in \cE(\lambda),\ u_{\mu}\in \cE(\mu)$ for distinct $\lambda,\mu \in C$. We claim that 
$$K_{v_{\lambda}u_{\mu}}=0\ \mbox{ if and only if }\ v_{\lambda}^t\, (x_1 \sigma_2-x_2 \sigma_1+x_0 \gamma) \, u_{\mu} \equiv 0
\mbox{ for all }x\in\PP^2.$$ 
Indeed, $K_{v_{\lambda}u_{\mu}}=0$ implies
$v_{\lambda}^t\, \sigma_1\, u_{\mu}=v_{\lambda}^t\, \sigma_2\, u_{\mu}=0$
by~(\ref{vin6.7}) and moreover $v_{\lambda}^t\, \gamma\, u_{\mu}=0$ by~(\ref{vin10}). The converse is obvious.

We mention few more implications of the definition of $K_{v_{\lambda}u_{\mu}}$:
\begin{itemize}
\item[i)]
If $\lambda,\mu$ are admissible points on $C$ then for every $v_{\lambda}\in \cE(\lambda)$ there exists exactly one  
$ v_{\mu}\in \cE(\mu)$ for which $K_{v_{\lambda}v_{\mu}}=0$.
\item[ii)]
For any $\lambda$ and $v_{\lambda}\in \cE(\lambda)$ there exists $\mu$ and $u_{\mu}\in \cE(\mu)$ such that 
$v_{\lambda}, u_{\mu}$ are admissible vectors.
\item[iii)]
Let $v_{\lambda}\in \cE(\lambda)$. Given any vector $u$ such that $K_{v_{\lambda} u}\neq 0$ there exist at most $d$ distinct points
$\mu^j\in C$ for which $u\in \cE(\mu^j)$ and $v_{\lambda}, u$ are admissible vectors.
\end{itemize}
Claim i) is obvious. For ii) recall that Span $\{\coker A (\mu);\, \mu\in C\} =k^{2 d}$. If $\mu$ didn't exist 
it would thus hold $v_{\lambda}^t\, \sigma_1\, u=v_{\lambda}^t\, \sigma_2\, u=0$ for all $u\in k^{2 d}$. This contradicts 
$\dim \{v_{\lambda}^t\, \sigma_1,\, v_{\lambda}^t\, \sigma_2\}^{\perp}=2d-2$. Claim iii) will follow from~(\ref{vin6.7}).
For any $t\in k$ 
$$\mu^j_i=\lambda_i-\frac{v_{\lambda}^t\, \sigma_i\, t u}{K_{v_{\lambda}\, t u}} ,\ i=1,2$$
we obtain the affine coordinates of $\mu^j$. But the defining polynomial of $C$ evaluated in $(1,\mu^j_1,\mu^j_2)$ 
has at most $d$ solutions for $\frac{t}{K_{v_{\lambda}\, t u}}$.

\begin{corollary}
For every $\lambda\in C$ there are at most finitely many $\mu$ that do not form an admissible pair of points.
\end{corollary}
 
\section{Elementary transformations act on the cokernel bundle of representation}
\label{secthree}

In this section we relate elementary transformations of vector bundles 
with the elementary transformations of pfaffian representations considered in
Section~\ref{beginsec}. 

\begin{definition} \label{defelembund} \rm{Let $\cE$ be a rank 2 vector bundle over a smooth and irreducible $C$. 
Take an effective reduced divisor
$Z$ on $C$ and consider the canonical surjection
$$\cE \rightarrow k(Z) \rightarrow  0,$$
where $k(Z)$ is a skyscraper sheaf at $Z$, i.e. rank 1 $\cO_{Z}$--module. Its kernel is a rank 2 vector bundle on $C$ called the 
\textit{elementary transformation of}
 $\cE$ \textit{at} $Z$. We denote it by $\cE'=\elem _Z(\cE)$.}
\end{definition}
\noindent
On $C$ it is equivalent to consider:
\begin{enumerate}
\item[(a)] Ruled surface $\pi: S \rightarrow C$ together with a base-point-free unisecant complete linear system $|H|$; 
\item[(b)] Rank 2 vector bundle $\cE$ over $C$ for which $S=\PP\cE$ and $\cE \cong \pi_{\ast}\cO_{\PP\cE}(H)$;
\item[(c)] Linearly normal scroll $R$ obtained as the image of the birational map $\phi_H: S \rightarrow R\subset \PP^N$
defined by $|H|$.
\end{enumerate}
Analogously we can define elementary transformation at a point $x\in C$ on each of the above:
\begin{enumerate}
\item[(a)] On the ruled surface $S$ we choose a point $s\in\pi^{-1}(x)$. Denote by $B$ the blow-up of $S$ at $s$. By Castelnuovo theorem we can contract 
the starting fibre $\pi^{-1}(x)$ in $B$ and obtain a new ruled surface $\pi': S'\rightarrow C$;
\item[(b)] $\cE'$ is obtained by Definition~\ref{defelembund} as the elementary transformation of $\cE$ at the divisor ${x}$ on $C$;
\item[(c)] Pick a point $r=\phi_H(s)$ on the scroll $R$ such that $\pi(s)=x$. Projection from $r$ yields a scroll $R'\subset \PP^{N-1}$.
\end{enumerate}
Fuentes and Pedreira~\cite{fuentes} checked that these definitions are compatible, namely $S'=\PP\cE'$. Moreover, $R'$ is the image
of $S'$ under the birational map defined by $|H'|$, where $H'=\nu^{\ast}(H)-\pi^{-1}(x)$. Here $\nu: S'\rightarrow S$ denotes the inverse of the
elementary transformation of $S$ at $s$. In particular $\cE' \cong \pi'_{\ast}\cO_{\PP\cE'}(H')$.

The inverse of an elementary transformation of a ruled surface is again an elementary transformation: if $S'$ is elementary transformation of 
$S$ at $s$, then $S$ is elementary transformation of 
$S'$ at $s'$, where $s'$ is the contraction of \\
(the exceptional divisor of the blow-up) $\cap$ (the original fibre $\pi^{-1}(x)$) $\in B$.\\
Up to tensoring line bundles, we can view the elementary transformation of a vector bundle
in Definition~\ref{defelembund} at $Z=x^1+\cdots+ x^m$ as an elementary transformation of $\PP \cE$ at $m$ points $s_i\in\pi^{-1}(x^i)$.
More precisely, there exists a skyscraper sheaf $k(Z)'$ that fits into the commutative diagram
\begin{equation} \label{diagelemtrvb}
\begin{array}{cl}
\cE\otimes \cO_{C}(-Z) & \\
 \downarrow  \! \substack{g} & \\
\cE' &  \! \! \! \!\! \! \! \!\! \! \! \!\! \! \! \!\! \!  \xrightarrow{e} \cE \rightarrow k(Z)\ .\\
\downarrow &\\
k(Z)' &
\end{array}\end{equation}
Here $\cE\otimes \cO_{C}(-Z)$ is the following elementary transformation of $\cE'$ at $Z$: Take a point $x\in C$ 
and free basis $\{X_1',X_2'\},\ \{X_1,X_2\}$ for $\cE'(x),\ \cE(x)$ respectively. Denote by 
$Y_i=e(X_i')=e_{1i} X_1+e_{2i}X_2,\ i=1,2$
the image of the morphism $e$. At every $x^i\in C,\ i=1,\ldots,m$ the determinant of $e$ is $0$, thus projectively $Y_1=Y_2=s_i\in\pi^{-1}(x^i)$ are
the points of the blow up of $\PP \cE$. 
Locally $k(Z)= k X_1 \oplus k X_2\, / \, kY_1+kY_2$, therefore $k(Z)$ is $0$ whenever $e$ is invertible and is a rank 1 bundle over the points $x^i$.  
Analogously $g$, given by the $2\times 2$ adjoint matrix of $e$, defines an elementary transformation of $\cE'$. 
It is easy to check that it equals $\cE\otimes \cO_{C}(-Z)$. For  $x^i\in C,\ i=1,\ldots,m$
the points of the blow up of $\PP \cE'$ equal $g_{11} X_1'+g_{21}X_2'=g_{12} X_1'+g_{22}X_2'$. Observe that 
$e(g_{11} X_1'+g_{21}X_2')=0$.
On the level of ruled surfaces $e^{-1}$ and $g^{-1}$ induce elementary transformations of
$\PP (\cE\otimes \cO_{C}(-Z))=\PP \cE$ and $\PP \cE'$ that are inverse to each other.\\


In the sequel we establish the explicit relation between elementary transformations of pfaffian representations and 
elementary transformations of the corresponding cokernel vector bundles. We use the notation of Definition~\ref{defelemrepr}.
Again we work in projective coordinates $x=(x_0,x_1,x_2)$ in $\PP^2$. 
\begin{theorem}\label{thmelemIIImur} 
Let $C$ be defined by $F=\pf (x_1 \sigma_2-x_2 \sigma_1+x_0 \gamma)$ and let 
$x_1 \sigma_2-x_2 \sigma_1+x_0 \tilde{\gamma},\ x_1 \sigma_2-x_2 \sigma_1+x_0 \bar{\gamma}$
be elementary transformations of Type I and II respectively. Denote by $\cE(x), \tilde{\cE}(x), \bar{\cE}(x)$ the corresponding cokernels
that are rank 2 vector bundles over $C$. The relating morphisms in the below diagrams can be expressed by 
elementary transformations of vector bundles,
$$\begin{array}{c}
\tilde{\cE} \\
S \uparrow\, R^t \downarrow    \\
\cE' \\
P \uparrow\ T^t \downarrow    \\
\cE
\end{array}
\ \ \ \ \ \mbox{and} \ \ \ \ \
\begin{array}{c}
 \phantom{x} \\
\left.   \begin{array}{c}
\phantom{x}  \\
\ \ Q \!\! \\
\phantom{x}  
\end{array}  \right\uparrow \\
\phantom{x}
\end{array}
\begin{array}{c}
\bar{\cE} \\
\downarrow    \\
\cE'' \\
\uparrow    \\
 \cE
\end{array}.
$$
\end{theorem}

\begin{proof}
For any parameters $t_1,t_2$ let $T(x),S(x)$ and $P(x),R(x)$ and $Q(x)$  be matrices with rational elements
\begin{eqnarray*}
T(x)=\Id+ \frac{x_0}{K_{v_{\lambda}u_{\mu}}\left(t_1(x_1-\lambda_1 x_0)+t_2(x_2-\lambda_2 x_0)\right)}
(t_1 \sigma_1 +t_2 \sigma_2)u_{\mu}v_{\lambda}^t,\\ 
S(x)=\Id+ \frac{x_0}{K_{v_{\lambda}u_{\mu}}\left(t_1(x_1-\lambda_1 x_0)+t_2(x_2-\lambda_2 x_0)\right)} 
u_{\mu}v_{\lambda}^t(t_1 \sigma_1 +t_2 \sigma_2),\\ 
P(x)=\Id+ \frac{x_0}{K_{v_{\lambda}u_{\mu}}\left(t_1(x_1-\mu_1 x_0)+t_2(x_2-\mu_2 x_0)\right)}
v_{\lambda}u_{\mu}^t(t_1 \sigma_1 +t_2 \sigma_2),\\
R(x)=\Id+ \frac{x_0}{K_{v_{\lambda}u_{\mu}}\left(t_1(x_1-\mu_1 x_0)+t_2(x_2-\mu_2 x_0)\right)} 
(t_1 \sigma_1 +t_2 \sigma_2)v_{\lambda}u_{\mu}^t\ \\
\mbox{and }Q(x)=\Id+\frac{2 \rho x_0}{t_1(x_1-\lambda_1 x_0)+t_2(x_2-\lambda_2 x_0)}
v_{\lambda} v_{\lambda}^t(t_1 \sigma_1 +t_2 \sigma_2).
\end{eqnarray*}
We will prove that $S(x)P(x)$ maps $\cE(x)$ into $\tilde{\cE}(x)$ and $Q(x)$ maps $\cE(x)$ into $\bar{\cE}(x)$.

First consider Type I transformations.
Observe that for $x$ different from $\lambda$ and $\mu$~(\ref{vin6.7}) implies
\begin{eqnarray*}
P(x)T^t(x)=T^t(x)P(x)=\Id\ & \mbox{and} \\
S^t(x)R(x)=R(x)S^t(x)=\Id, & 
\end{eqnarray*}
Using this,~(\ref{vin14}) rewrites into
\begin{eqnarray*}
T(x)(x_1 \sigma_2-x_2 \sigma_1+x_0\gamma)T^t(x) &= &S^t(x)(x_1 \sigma_2-x_2 \sigma_1+x_0\tilde{\gamma})S(x), 
\end{eqnarray*}
or equivalently
\begin{equation} \label{vin2319}
R(x)T(x)(x_1 \sigma_2-x_2 \sigma_1+x_0\gamma)=(x_1 \sigma_2-x_2 \sigma_1+x_0\tilde{\gamma})S(x)P(x).
\end{equation}
Therefore
\begin{equation}
\label{vin23}
\tilde{\cE}(x)\,R(x)T(x) \subseteq \cE(x)\ \mbox{ and }\ S(x)P(x)\, \cE(x) \subseteq \tilde{\cE}(x).
\end{equation}
Transposing~(\ref{vin2319}) and multiplying by the inverse matrices yields the reverse inclusions.
This means that the rational matrix function
$S(x)P(x)$
maps the vector bundle corresponding to the starting representation
to the cokernel bundle of the new representation of Type I. 

Similarly, $\Id+\frac{2 \rho x_0}{t_1(x_1-\lambda_1 x_0)+t_2(x_2-\lambda_2 x_0)}
(t_1 \sigma_1 +t_2 \sigma_2)v_{\lambda} v_{\lambda}^t$ is the inverse of $Q^t(x)$ for $x\neq \lambda$. 
For Type II elementary transformations,~(\ref{vin14sec}) 
becomes
$$Q^t(x)^{-1}(x_1 \sigma_2-x_2 \sigma_1+x_0\gamma)=
(x_1 \sigma_2-x_2 \sigma_1+x_0 \bar{\gamma})Q(x)$$
and
\begin{equation}
\label{vin23a}
Q(x)\, \cE(x)=\bar{\cE}(x).
\end{equation}
Thus the rational matrix function $Q(x)$ maps $\cE(x)$ to the cokernel bundle $\bar{\cE}(x)$ of the representation 
obtained by the Type II transformation.

An attentive reader will notice that matrices $S(x),T(x),P(x),R(x),Q(x)$ also depend on parameters $t_1,t_2$. But for $x$ 
such that 
$t_1(x_1-\lambda_1 x_0)+t_2(x_2-\lambda_2 x_0) \neq 0,\ t_1(x_1-\mu_1 x_0)+t_2(x_2-\mu_2 x_0) \neq 0$ and for every 
$\varepsilon(x)\in\cE(x),\,\tilde{\varepsilon}(x)\in\tilde{\cE}(x),\,\bar{\varepsilon}(x)\in\bar{\cE}(x)$ the following vectors 
\begin{eqnarray*}
P(x)\varepsilon(x),\ \  \tilde{\varepsilon}^t(x)T(x),\ \  S(x)\varepsilon (x),\ \ \tilde{\varepsilon}^t(x)R(x)\ \mbox{ and }\ 
Q(x)\varepsilon(x),\ \  Q^{-1}(x)\bar{\varepsilon}(x)
\end{eqnarray*}
are independent of $t_1,t_2$. This is an immediate corollary of~(\ref{vin9a}), (\ref{vin13a}) and (\ref{vintype2})
since 
\begin{eqnarray*}
0&=&u_{\mu}^t(x_1 \sigma_2-x_2 \sigma_1+x_0\gamma)\varepsilon(x)=
u_{\mu}^t\left((x_1-\mu_1 x_0) \sigma_2-(x_2-\mu_2 x_0) \sigma_1\right)\varepsilon(x),\\
0&=&\tilde{\varepsilon}^t(x)(x_1 \sigma_2-x_2 \sigma_1+x_0\tilde{\gamma})u_{\mu}=
\tilde{\varepsilon}^t(x)\left((x_1-\lambda_1 x_0) \sigma_2-(x_2-\lambda_2 x_0) \sigma_1\right)u_{\mu},\\
0&=&v_{\lambda}^t(x_1 \sigma_2-x_2 \sigma_1+x_0\gamma)\varepsilon(x)=
v_{\lambda}^t \left((x_1-\lambda_1 x_0) \sigma_2-(x_2-\lambda_2 x_0) \sigma_1\right)\varepsilon(x),\\
0&=&\tilde{\varepsilon}^t(x)(x_1 \sigma_2-x_2 \sigma_1+x_0\tilde{\gamma})v_{\lambda}=
\tilde{\varepsilon}^t(x) \left((x_1-\mu_1 x_0) \sigma_2-(x_2-\mu_2 x_0) \sigma_1\right)v_{\lambda}\ \\
\mbox{and } 0&=&v_{\lambda}^t(x_1 \sigma_2-x_2 \sigma_1+x_0\bar{\gamma})\bar{\varepsilon}(x)=
v_{\lambda}^t\left((x_1-\lambda_1 x_0) \sigma_2-(x_2-\lambda_2 x_0) \sigma_1\right)\bar{\varepsilon}(x).
\end{eqnarray*}
Hence $P(x),R(x),\ S(x),T(x)$ and $Q(x)$ do not depend on $t_1,t_2$ when restricted to $\cE(x),\ \tilde{\cE}(x)$ 
and $\bar{\cE}(x)$ respectively.
Moreover, by using~(\ref{vin6.7}) for the admissible pair of vectors
$v_{\lambda}\in\cE(\lambda)\cap\tilde{\cE}(\mu),\,u_{\mu}\in\cE(\mu)\cap\tilde{\cE}(\lambda)$,
we get 
\begin{equation}\label{lammuzeros}
 P(\lambda)v_{\lambda}=0, \ \ \  v_{\lambda}^t T(\mu)=0,\ \ \ 
u_{\mu}^t R(\lambda) =0,\ \ \ S(\mu) u_{\mu}=0.
\end{equation}

Next we analyse how elementary transformations act on rational sections of the cokernel bundles.
Let $q(x)$ be a \textit{rational section} of $\cE(x)$. This is a $2d-$tuple of rational functions along $C$ with the property $q(x)\in \cE(x)$. 
For example, by~(\ref{pfaflinalg}) it can be obtained by dividing a column in the adjoint matrix of the representation by a degree $d-1$ polynomial. 
Then by~(\ref{vin23}) and~(\ref{vin23a}) 
\begin{eqnarray}\label{vin31}
\tilde{q}(x)=S(x)P(x)\, q(x)\ \mbox{ and }\ \bar{q}(x)=Q(x)\, q(x)
\end{eqnarray}
are rational sections of $\tilde{\cE}(x)$ and  $\bar{\cE}(x)$ respectively. 
Additionally, denote by $q'(x)=P(x)\, q(x)=R^t(x)\, \tilde{q}(x)$ sections in the auxiliary bundle $\cE'(x)=P(x)\, \cE(x)=R^t(x)\, \tilde{\cE}(x)$.
Proposition~\ref{propratsect} will discuss the order of these rational sections at various points on $C$. This will describe 
$P,S$ and $Q$ in terms of elementary transformations of vector bundles considered in~(\ref{diagelemtrvb}). 
In other words, 
\begin{itemize}
\item
$P$ at $\lambda$ equals the morphism in the diagram
$$ \begin{array}{cl}
\cE'\otimes \cO_{C}(-\lambda) & \\
 \downarrow  & \\
\cE &  \! \! \! \!\! \! \! \!\! \! \! \!\! \! \! \!\! \!  \xrightarrow{P} \cE' \rightarrow k(\lambda).
\end{array}$$
$\cE'\otimes \cO_{C}(-\lambda)$ (the induced elementary transformation of $\cE$ at $\lambda$) 
is by (\ref{lammuzeros}) obtained through the blow--up of $v_{\lambda}\in \PP\cE(\lambda)$;
\item
$P$ at $\mu$ equals the elementary transformation of $\cE$ at $\mu$ in the diagram
$$ \begin{array}{cl}
\cE\otimes \cO_{C}(-\mu) & \\
 \downarrow   & \\
\cE' &  \! \! \! \!\! \! \! \!\! \! \! \!\! \! \! \!\! \!  \xrightarrow{T^t} \cE \rightarrow k(\mu).
\end{array}$$
Again by (\ref{lammuzeros}), $\cE\otimes \cO_C(-\mu)$ is obtained through the blow--up of $v_{\lambda}\in \PP\cE'(\mu)$.

\end{itemize}
The same way  $S$ consists of the following, 
\begin{itemize}
\item
$S$ at $\lambda$ equals the elementary transformation of ${\cE'}$ at $\lambda$ in the diagram
$$ \begin{array}{cl}
{\cE'}\otimes \cO_{C}(-\lambda) & \\
 \downarrow   & \\
\tilde{\cE}  &  \! \! \! \!\! \! \! \!\! \! \! \!\! \! \! \!\! \!  \xrightarrow{R^t} {\cE'} \rightarrow k(\lambda)^{\sim},
\end{array}$$
where ${\cE'}\otimes \cO_C(-\lambda)$ is obtained through the blow--up of $u_{\mu}\in \PP \tilde{\cE}(\lambda)$;
\item
$S$ at $\mu$ equals the morphism in the diagram
$$ \begin{array}{cl}
\tilde{\cE} \otimes \cO_{C}(-\mu) & \\
 \downarrow   & \\
{\cE'} &  \! \! \! \!\! \! \! \!\! \! \! \!\! \! \! \!\! \!  \xrightarrow{S} \tilde{\cE}  \rightarrow k(\mu)^{\sim},
\end{array}$$
where $\tilde{\cE} \otimes \cO_C(-\mu)$ is obtained through the blow--up of $u_{\mu}\in \PP {\cE'}(\mu)$.
\end{itemize}
On the other hand, $Q$ and $Q^{-1}$ are described by the following diagram of elementary transformations at $\lambda \in \pic C$
$$ \begin{array}{ccccccc}
     & & \bar{\cE} \otimes \cO_{C}(-\lambda) \! \! \! \! \! \!&                          &\! \! \! \! \! \! {\cE} \otimes \cO_{C}(-\lambda) & & \\
                    &  &           &  \searrow  \ \ \ \swarrow                                        &  & & \\
& \ \ \ \ \ \ \ \ \ \ \ \ \ \ \ \ \ \ \ \ \ \ k(\lambda)^{-} \longleftarrow \! \! \! \! \! \! \! \! \! \! \! \! \! \! \! \! \! \! \! \! \! \!
& \ \ \ \ \ \ \ \ \ \ \ \bar{\cE}  \longleftarrow 
& {\cE''} &   \longrightarrow  {\cE} \ \ \ \ \ \ \ \ \ \ \ &\! \! \! \! \! \! \! \! \! \! \! \! \! \! \! \! \! \! \! \! \! \!
\longrightarrow k(\lambda), \ \ \ \ \ \ \ \ \ \ \ \ \ \ \ \ \ \ \ \ \ \ & 
\end{array}$$
where $\cE''$ is an elementary transformation of both $\cE$ and $\bar{\cE}$.
This will conclude the proof of Theorem~\ref{thmelemIIImur}.
\end{proof}

\begin{figure}[ht] 
\begin{center}
\includegraphics*[100pt, 15pt][10cm, 9cm]{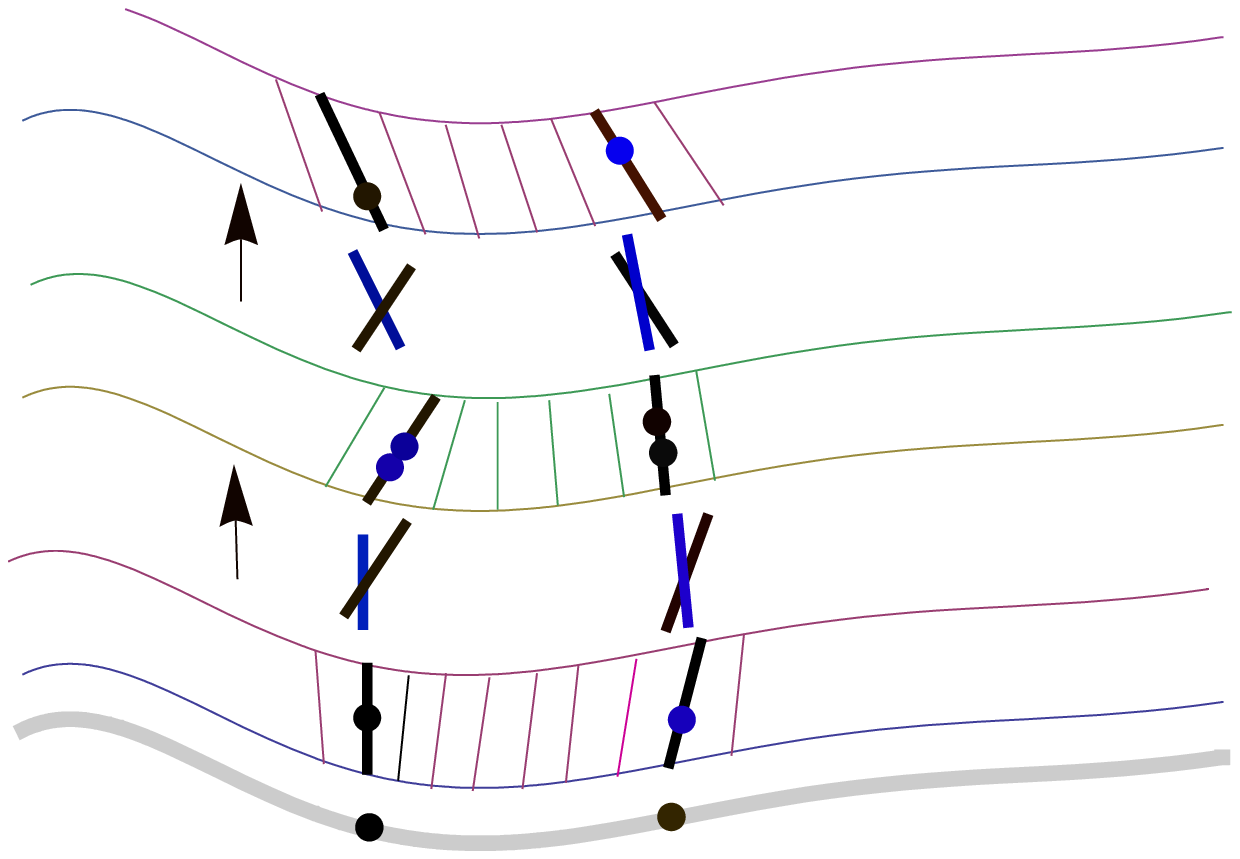}
  \put(-44,1){${\mu}$}
  \put(-144,-3){${\lambda}$}
  \put(-20,44){$\cE$}
  \put(-48,124){$u_{\mu}$}
  \put(-148,46){$v_{\lambda}$}
  \put(-20,124){$\cE'$}
  \put(-20,204){$\tilde{\cE}$}
  \put(-46,114){$v_{\lambda}$}
  \put(-152,196){$u_{\mu}$}
  \put(-170,80){$P$}
  \put(-170,160){$S$}
  \put(-20,6){$C$}
\hspace{2cm}
\includegraphics*[200pt, 15pt][9.5cm, 9cm]{SlikaEekTransf.eps}
  \put(-38,-4){${\lambda}$}
  \put(-10,2){$C$}
  \put(-15,200){$\ \bar{\cE}$}
  \put(-15,122){${\ \cE''}$}
  \put(-12,46){${\ \cE}$}
 \put(-73,78){$Q$}
  \put(-81,140){$\Uparrow$}
  \put(-80,130){$\Arrowvert$}
  \put(-80,120){$\Arrowvert$}
  \put(-80,110){$\Arrowvert$}
  \put(-80,100){$\Arrowvert$}
  \put(-80,90){$\Arrowvert$}
  \put(-80,80){$\Arrowvert$}
\label{pict6}
\end{center}
\caption{Elementary transformations of Type I and II}
\end{figure}

\begin{proposition} 
\label{propratsect}
Let $q(x)$ be a rational section of $\cE(x)$ and $o_{\xi}(q)$ its order at $\xi$ and let $\tilde{q}(x),\,  \bar{q}(x),\,  q'(x)$ be sections 
of $\tilde{\cE}(x),\, \bar{\cE}(x),\, \cE'(x)$ respectively.

If $\xi$ is different from $\lambda$ and $\mu$ then $o_{\xi}(q)=o_{\xi}(q')=o_{\xi}(\tilde{q})=o_{\xi}(\bar{q})$. 

If $\xi=\lambda$ then: \\
$o_{\lambda}(q')=o_{\lambda}(q) +1$ iff $q$ at $\lambda$ passes through $v_{\lambda}$, else 
$o_{\lambda}(q)=o_{\lambda}(q')$;\\
$o_{\lambda}(q')=o_{\lambda}(\tilde{q}) +1$ iff $\tilde{q}$ at $\lambda$ passes through $u_{\mu}$, else
$o_{\lambda}(\tilde{q})=o_{\lambda}(q')$;\\
$o_{\lambda}(q)=o_{\lambda}(\bar{q}) $ iff $q$ and $\bar{q}$ at $\lambda$ both pass through $v_{\lambda}$, 
else either $o_{\lambda}(q)=o_{\lambda}(\bar{q})+1$ or $o_{\lambda}(\bar{q})=o_{\lambda}(q)+1$ .

If $\xi=\mu$ then:\\
$o_{\mu}(q)=o_{\mu}(q')+1$ iff $q'$ at $\mu$ passes through $v_{\lambda}$, else 
$o_{\mu}(q)=o_{\mu}(q')$;\\
$o_{\mu}(\tilde{q})=o_{\mu}(q') +1$ iff $q'$ at $\mu$ passes through $u_{\mu}$, else $o_{\mu}(\tilde{q})=o_{\mu}(q')$.
\end{proposition}

\begin{proof}
Pick a local parameter $t$ at $\xi=(\xi_0,\xi_1,\xi_2)\in C$ and expand $x$ and $q(x)$
into the power series with center $\xi$. We can assume that $\xi_0=1$ and $\frac{\partial F}{\partial x_0}(\xi)=0$. Thus
\begin{eqnarray*}
x(t)&=&\left( 1+\cdots,\xi_1+t\frac{\partial F}{\partial x_1}(\xi)+\cdots,
\xi_2+t\frac{\partial F}{\partial x_2}(\xi)+\cdots \right), \\
q(x)&=&c_o t^o+c_{o+1}t^{o+1}+\cdots,
\end{eqnarray*}
where $o=o_{\xi}(q)$ and $\cdots$ denote higher order terms. By definition $ c_o\neq 0$.  Since
$(x_1 \sigma_2-x_2 \sigma_1+x_0\gamma)\,q(x)= 0$ for all $x\in C$, in terms of the local parameter $t$ we get 
$t^o\, ((\xi_1+\cdots) \sigma_2-(\xi_2+\cdots) \sigma_1+(\xi_0+\cdots)\gamma)\, (c_o +c_{o+1}t+\cdots)\equiv 0$ which implies
$$c_o\in \cE(\xi).$$
Analogously we expand
$${q'}(x)=c_{{o'}} t^{{o'}}+\cdots,\ \ \ \tilde{q}(x)=c_{\tilde{o}} t^{\tilde{o}}+\cdots,\ \ \ \bar{q}(x)=c_{\bar{o}} t^{\bar{o}}+\cdots,$$
where ${o'}=o_{\xi}({q'}),\ \tilde{o}=o_{\xi}(\tilde{q}),\ \bar{o}=o_{\xi}(\bar{q})$
and $c_{\tilde{o}}\in \tilde{\cE}(\xi),\ c_{\bar{o}}\in \bar{\cE}(\xi).$

Let $\xi\neq\lambda,\mu$. Then $S(x),P(x)$ and $Q(x)$ are defined at $\xi$ and can be expanded along $t$. 
Thus~(\ref{vin31}) gives
$${q'}(x)= P(\xi)\, c_o\, t^o+\cdots, \ \tilde{q}(x)= S(\xi)P(\xi)\, c_o\, t^o+\cdots \ \mbox{ and }\ \bar{q}(x)= Q(\xi)\, c_o\, t^o+\cdots  .$$
This proves $o'\geq o,\ \tilde{o}\geq o$ and $\bar{o}\geq o$. 
The power series for $q'(x),\ \tilde{q}(x),\ \bar{q}(x)$ and   
$q(x)=T^t(x)\, {q'}(x),\ q(x)=T^t R^t(x)(x)\, \tilde{q}(x),\ q(x)=Q(x)^{-1} \bar{q}(x)$ analogously imply 
$o\geq o',\ o\geq \tilde{o},\ o\geq \bar{o}$.

Let now $\xi=\lambda$. $P(x)$ and $R(x)$ are defined at $\lambda$ since $\lambda\neq \mu$. 
By the considerations in the proof of~(\ref{lammuzeros}) we notice that 
$P(\lambda)\, c_o=0$ if and only if $c_o$ is a multiple of $v_{\lambda}$ and that
$R^t(\lambda)\, c_{\tilde{o}}=0$ if and only if $c_{\tilde{o}}$ is a multiple of $u_{\mu}$.
Like above, ${q'}(x)= P(\lambda)\, c_o\, t^o+\cdots$ implies $o'\geq o$. Moreover, strict inequality holds iff 
$c_o$ is a multiple of $v_{\lambda}$ (i.e., $q$ at $\lambda$ passes through $v_{\lambda}$). 
The same way ${q'}(x)= R^t(\lambda)\, c_{\tilde{o}}\, t^{\tilde{o}}+\cdots$ implies $o'\geq \tilde{o}$ and strict inequality holds iff
$\tilde{q}$ at $\lambda$ passes through $u_{\mu}$. \\
Recall that $\lambda=(1,\lambda_1,\lambda_2)$ is a regular point on 
$C$, therefore one of $\frac{\partial F}{\partial x_1}(\lambda),\ \frac{\partial F}{\partial x_2}(\lambda)$ must be nonzero.
Then
\begin{eqnarray*}
T(x)&=&\Id+ \frac{1}{t\, K_{v_{\lambda}u_{\mu}}
\left(t_1(\frac{\partial F}{\partial x_1}(\lambda)+\cdots )+t_2(\frac{\partial F}{\partial x_2}(\lambda)+\cdots )\right)} 
(t_1 \sigma_1 +t_2 \sigma_2)u_{\mu}v_{\lambda}^t,\\ 
S(x)&=&\Id+ \frac{1}{t\, K_{v_{\lambda}u_{\mu}}
\left(t_1(\frac{\partial F}{\partial x_1}(\lambda)+\cdots )+t_2(\frac{\partial F}{\partial x_2}(\lambda)+\cdots )\right)} 
u_{\mu}v_{\lambda}^t(t_1 \sigma_1 +t_2 \sigma_2) 
\end{eqnarray*}
imply
\begin{eqnarray*}
q(x)=T^t(x)\, q'(x)\!\!&\!\!=\!\!&\! \!\frac{-1}{ K_{v_{\lambda}u_{\mu}}
\left(t_1 \frac{\partial F}{\partial x_1}(\lambda)+t_2 \frac{\partial F}{\partial x_2}(\lambda)\right)}\, 
v_{\lambda}u_{\mu}^t\, (t_1 \sigma_1 +t_2 \sigma_2)\, c_{o'}\ t^{o'-1}+\cdots \! ,\\
\tilde{q}(x)=S(x)\, q'(x)\!\!&\!\!=\!\!&\!\!\frac{1}{ K_{v_{\lambda}u_{\mu}}
\left(t_1 \frac{\partial F}{\partial x_1}(\lambda)+t_2 \frac{\partial F}{\partial x_2}(\lambda)\right)}\, 
u_{\mu}v_{\lambda}^t\,(t_1 \sigma_1 +t_2 \sigma_2)\, c_{o'}\ t^{o'-1}+\cdots \! .
\end{eqnarray*}
This proves $o\geq o'-1$ and $\tilde{o}\geq o'-1$. \\
The same arguments apply to
\begin{eqnarray*}
\bar{q}(x)=Q(x)\, q(x)=\frac{2 \rho }{ t_1\frac{\partial F}{\partial x_1}(\lambda)+t_2 \frac{\partial F}{\partial x_2}(\lambda) }\,
v_{\lambda} v_{\lambda}^t (t_1 \sigma_1 +t_2 \sigma_2)\, c_{{o}} \ t^{{o}-1}+\cdots , \\
q(x)=Q(x)^{-1}\, \bar{q}(x)=\frac{-2 \rho }{ t_1\frac{\partial F}{\partial x_1}(\lambda)+t_2 \frac{\partial F}{\partial x_2}(\lambda) }\,
v_{\lambda} v_{\lambda}^t (t_1 \sigma_1 +t_2 \sigma_2)\, c_{\bar{o}} \ t^{\bar{o}-1}+\cdots  .
\end{eqnarray*}
This proves $\bar{o} \geq o-1$ and $o\geq \bar{o}-1$. Observe also that at least one of $q,\bar{q}$ 
passes through $v_{\lambda}$, which is equivalent to $c_{{o}}$ or $c_{\bar{o}}$ being a multiple of $v_{\lambda}$.  
More precisely, $o=\bar{o}$ iff $q$ and $\bar{q}$ at $\lambda$ both pass through $v_{\lambda}$, $\bar{o}=o-1$ iff 
 $c_o\neq v_{\lambda}\in \PP\cE$ and $o=\bar{o}-1$ iff $c_{\bar{o}}\neq v_{\lambda}\in \PP\bar{\cE}$. 

It remains to consider the point $\xi=\mu$. 
Since $T^t(x)$ and $S(x)$ are defined at $\mu$, $q(x)=T^t(x)\, q'(x)$ and $\tilde{q}(x)=S(x)\, q'(x)$ imply $o\geq o'$ and 
$\tilde{o}\geq o'$. Moreover, by~(\ref{lammuzeros}) $o> o'$ iff $q'$ passes through $v_{\lambda}$ and  
$\tilde{o}> o'$ iff $q'$ passes through $u_{\mu}$. Since $P(x)$ and $R^t(x)$ have poles at $\mu$, $q'(x)=P(x)\, q(x)=R^t(x)\, \tilde{q}(x)$
imply $o'\geq o-1$ and $o'\geq \tilde{o}-1$.
\end{proof}

\section{Pfaffians arising from decomposable vector bundles}
\label{decsec}

The \textit{Kummer variety} $\cK_C$ of $C$ is by definition the quotient of the Jacobian $J C$ by the involution 
$\cL \mapsto \cL^{-1}\otimes \cO_C(d-3)$. It is standard~\cite{laszlo}, that $M_C(2,\cO_C(d-3))$ is singular along the Kummer variety.

In~\cite[Section 5]{anitapfaff} we established a one to one correspondence between decomposable vector bundles in 
$M_C(2,\cO_C(d-3))\ \backslash \ \Theta_{2,\cO_C(d-3)}$ and the open subset of Kummer variety
$\cK_C\ \backslash \ W$, 
where $W=\{$line bundles of degree $\frac{1}{2}d(d-3)$ with no sections$\}$. 
On the other hand, the cokernel $\cE$ of a pfaffian representation is decomposable if and only if it is of the form
$$\cE \cong \coker M \oplus \left((\coker M)^{-1}\otimes \cO_C(d-1)\right)$$
for $M$ a determinantal representation of $C$.
The line bundle $\coker M$ is of degree $\frac{1}{2}d(d-1)$ with $H^0(C,\coker M\otimes \cO_C(-1))=0$ and 
$\cE$ is the cokernel of the \textit{decomposable matrix}
$$\left[ \begin{array}{cc}
0 & M \\
-M^t & 0
\end{array}\right].$$
Denote $\cL=\coker M \otimes \cO_C(-1)$  and $\cL_t=\coker M^t \otimes \cO_C(-1)$
that both have degree $\frac{1}{2}d(d-3)$. Then
$$\cL \oplus \cL_t \ \in \ M_C(2,\cO_C(d-3))\ \backslash \ \Theta_{2,\cO_C(d-3)}$$
and by the above
$$\cL_t \ \cong \  \cL^{-1}\otimes \cO_C(d-3) ,$$
which explains the involution on the Kummer variety.
The involution appears, because in general $M$ and $M^t$ are not equivalent determinantal representations, but 
$$\left[ \begin{array}{cc}
0 & M \\
- M^t & 0
\end{array}\right]\ \mbox{ and }\ \left[ \begin{array}{cc}
0 & M^t \\
- M & 0
\end{array}\right]=\left[ \begin{array}{cc}
0 & \Id \\
- \Id & 0
\end{array}\right]\, \left[ \begin{array}{cc}
0 & M \\
- M^t & 0
\end{array}\right]\, \left[ \begin{array}{cc}
0 & -\Id \\
\Id & 0
\end{array}\right]$$
are equivalent pfaffian representations.

The proof of Theorem~\ref{thmelemIIImur} and~\cite[Theorem 4]{vinnikovElTr} imply the following
\begin{corollary}
\label{cordec}
Let $\cE$ be the decomposable cokernel of a decomposable pfaffian representation and
let $D=\div \coker M$ and $D'=\div \coker M^t$. Pick an admissible array of the form
$$\left[ \begin{array}{c}
v\\
0
\end{array}\right]
\in\cE(\lambda),\ 
\left[ \begin{array}{c}
0\\
u
\end{array}\right]\in\cE(\mu)$$ 
and perform the Type I elementary transformation.
Then 
$$\begin{array}{c}
\cE \cong \cO_C(D) \oplus \cO_C(D')\\
P \downarrow \\
\cE' \cong \cO_C(D+\lambda-\mu) \oplus \cO_C(D')\\
S \downarrow \\
\tilde{\cE} \cong \cO_C(D+\lambda-\mu) \oplus \cO_C(D'+\mu-\lambda).
\end{array}$$
\end{corollary}

\begin{example}\label{remkdec}\rm{
Recall the explicit calculation of pfaffian representations in the canonical form~(\ref{seccanform}) in Remark~\ref{decremex}.
Using local parameters and the implicit function theorem, it is easy to see that by
setting the entries $a^0_{ij}$ for $ 1\leq i<j\leq d$ and $d\leq i<j\leq 2d$ in $A_0$ to zero,
represents the singular locus of the space of all pfaffian representations of $C$. 
Note that these representations are decomposable and non-equivalent 
to each other.
Moreover, Vinnikov's explicit parametrisation of determinantal representations of $C$ proves that these are 
all the decomposable pfaffian representations.
}\end{example}
We will conclude this section by
\begin{theorem}\label{cortwoparts}
From any given pfaffian representation of $C$ we can build all the nonequivalent pfaffian representations of $C$ by 
finite sequences of Type I and Type II elementary transformations.
\end{theorem}

\begin{proof}
A finite sequence of Type I and Type II elementary transformations gives us a new pfaffian representation of $C$ that 
is in general not equivalent to the starting one. This follows from Theorem~\ref{thmelemIIImur} since the cokernel bundles 
$\cE(x),\ \tilde{\cE}(x),\, \bar{\cE}(x)$ are in general not isomorphic. On the other hand, the auxiliary bundles 
$\cE'(x)$ and $\cE''(x)$ have determinants different to $\cO_C(d-1)$ and thus can not be cokernels of pfaffian representations.

Pick the cokernel bundle $\cE$ of $x_1 \sigma_2-x_2 \sigma_1+x_0 \gamma$. 
We will assume that the representation is in the canonical form~(\ref{seccanform}).

In the first step we bridge $\cE$ with a decomposable vector bundle by applying a finite number of Type II elementary transformations.
A finite sequence of $m$ Type II elementary transformations by recursion yields a new representation
$x_1 \sigma_2-x_2 \sigma_1+x_0 \gamma_{m}$, where
$$\gamma_m=\gamma+ \sum_{j=1}^m \rho_j \ \sigma_2 v_{\lambda^j} \wedge \sigma_1 v_{\lambda^j} .$$
The above constants $\rho_j\in k$ and points $\lambda^j\in C$ are arbitrary and 
$v_{\lambda^j}\in \cE_{j-1}(\lambda^j):=\coker (\lambda^j_1 \sigma_2-\lambda^j_2 \sigma_1+\lambda^j_0 \gamma_{j-1})$
with $\gamma_0=\gamma$.
Since every union $\{\cE_{j}(x)\}_{x\in C}$ 
spans the whole $k^{2d}$, we can (by suitable choices of $v_{\lambda^j}$) generate enough independent rank 2 
matrices $\sigma_2 v_{\lambda^j} \wedge \sigma_1 v_{\lambda^j}$ whose linear combination will yield a decomposable matrix
$\gamma_m$. Indeed, 
$$\sigma_2 v_{\lambda^j} \wedge \sigma_1 v_{\lambda^j}=
\left[\! \! \! \! \! \!  \begin{array}{cc}
\begin{array}{ccc} 0 & & \\ & \ddots &  \! \! \! \! \! \! \! \phantom{x}^{v_{d+n} v_{d+l} (p_l-p_n)}\\  
& &  \! \! \! \! \! \! \! \! \! \! \! \! \! \! \! \! \! \! 0\end{array} & \ast \\
-\ast^t & \begin{array}{ccc} 0 & & \\ & \ddots & \! \! \! \! \! \! \!  \phantom{x}^{v_{n} v_{l} (p_l-p_n)} \\  
& & \! \! \! \! \! \! \! \! \! \! \! \! \! \! \! \! \! \! 0\end{array}
\end{array}\! \! \! \! \! \! \! \right]$$
form a basis for the two diagonal square blocks of $2d\times 2d$ skew--symmetric matrices.  
Here $v_i, i=1,\ldots,2d$ are the components of $v_{\lambda^j}$ and $p_i, i=1,\ldots,d$ are the 
distinct entries in the diagonal matrix $D$ in $\sigma_1$.

In the second step we bridge any two decomposable cokernel bundles $_d\dot{\cE},\ _d\ddot{\cE}$ 
by applying a finite number of Type I elementary transformations.
We can write 
\begin{eqnarray*}
_d\dot{\cE}=\cO_C(\dot{D})\oplus \cO_C(H-\dot{D}),\\
_d\ddot{\cE}=\cO_C(\ddot{D})\oplus \cO_C(H-\ddot{D}),\end{eqnarray*}
where $H$ is the divisor of a degree $d-1$ polynomial in $\PP^2$.
Let $g$ be the genus of $C$. For general points $\lambda^1,\ldots,\lambda^g\in C$ the Riemann--Roch theorem implies
$$\dot{D}-\ddot{D}+\lambda^1+\cdots +\lambda^g \equiv \mu^1+\cdots +\mu^g$$
for some distinct points $\mu^1,\ldots,\mu^g \in C$. Vinnikov~\cite[Theorems 5 and 6]{vinnikovElTr} proved that 
the indices of $\mu^i$'s can be permuted in such a way that for $n=1,\ldots,g$ every
$$_d\cE_{n}:=\cO_C\left( \dot{D}+ \sum_{i=1}^n \lambda^i-\mu^i\right) \oplus \cO_C \left(H-\dot{D}+\sum_{i=1}^n \mu^i-\lambda^i \right)$$
is the cokernel bundle of a decomposable pfaffian representation of $C$. By recursion and Corollary~\ref{cordec}, 
$_d\cE_{n}$ is obtained from $_d\cE_{n-1}$ by the Type I elementary transformation based on
$$\left[ \begin{array}{c}
v_n\\
0
\end{array}\right]
\in \, _d\cE_{{n-1}}(\lambda^n),\ 
\left[ \begin{array}{c}
0\\
u_n
\end{array}\right]\in \, _d\cE_{{n-1}}(\mu^n),$$
where we put $_d\cE_{0}=\, _d\dot{\cE}$.

Recall that the inverses of elementary transformations are again elementary transformations of the same type. This concludes the proof.
\end{proof}

\section{Plane quartic}
\label{seclast}

A nonsingular plane quartic $C$ is a non-hyperelliptic genus 3 curve embedded by its
canonical linear system $|K_C|$.
We parametrised 
$M_C(2,\cO_C(1))\ \backslash \ \Theta_{2,\cO_C(1)}$
by pfaffian representations of $C$. 
The moduli space
$M_C(2,\cO_C(1))\cong M_C(2,\cO_C)$
can be embeded as a Coble quartic hypersurface in $\PP^7$ with singularities along the 3--dimensional Kummer variety $\cK_C$. 
For references 
check~\cite{narasimhanram},~\cite{laszlo},~\cite{beauvilleTh}.

In this section we establish the connection between two distinguished kinds of rank 2 vector bundles on $C$, 
namely decomposable bundles corresponding to even theta characteristics and indecomposable Aronhold bundles.

A \textit{theta characteristic} of $C$ is a line bundle 
$\cL_{\vartheta}$ with the property 
$$\cL_{\vartheta}^{\otimes 2}\cong \omega_C \cong \cO_C(1).$$ 
A theta characteristic is called even (odd) if $\dim H^0(C,\cL_{\vartheta})$ is even (odd).
By~\cite{dolgachev} there are exactly 36 even theta characteristics on a smooth plane quartic. Since $C$ 
is not hyperelliptic no even theta characteristic vanishes, which means $H^0(C,\cL_{\vartheta})=0$. Therefore
$$H^0(C,\cL_{\vartheta}\oplus \cL_{\vartheta})=0\ \mbox{ and } \det \cL_{\vartheta}\oplus \cL_{\vartheta}\cong \omega_C.$$
In our notation
$$ \cL_{\vartheta}\oplus \cL_{\vartheta} \ \in \ M_C(2,\cO_C(1))\ \backslash \ \Theta_{2,\cO_C(1)}$$
and thus every even theta characteristic induces a decomposable pfaffian representation of $C$.
By Section~\ref{decsec} the corresponding pfaffian representations are 
$$\left[ \begin{array}{cc}
0 & M_{\vartheta} \\
- M_{\vartheta} & 0
\end{array}\right],$$
where $\cL_{\vartheta}=\coker M_{\vartheta} \otimes \cO_C(-1)$ and $M_{\vartheta}=M_{\vartheta}^t$.

\begin{example} \rm{
An easy computation in \texttt{Wolfram Mathematica} shows that, if we  
add $a^0_{i\, d+j}=a^0_{j\, d+i}$ for $ 1\leq i<j\leq d$ 
to the equations describing the decomposable representations of Example~\ref{remkdec} 
in the canonical form~(\ref{seccanform}) ,
we get exactly 36 solutions. As expected, they are the 36 $M_{\vartheta}$.}
\end{example}

These considerations generalise to the following proposition.
\begin{proposition} \label{thetasymmdet}
For a line bundle $\cL$ on a nonsingular plane quartic $C$ with $H^0(C,\cL)=0$ the following are equivalent:
\begin{itemize}
\item $\cL$ is an even theta characteristics on $C$,
\item $\cL \ \cong \  \cL^{-1}\otimes \cO_C(1)$,
\item $\cL=\coker M \otimes \cO_C(-1)$ where $M$ is a symmetric determinantal representation 
of $C$ with the property $\coker M\cong \coker M^t$.
\end{itemize}
\end{proposition}  
  
To every symmetric determinantal representation $M_{\vartheta}$ of $C$ one can associate a net of quadrics $\mathbb{M}_{\vartheta}$ in $\PP^3$.
The base locus of $\mathbb{M}_{\vartheta}$ consists of 8 points $b_1,\ldots,b_8$, called the Cayley octad. We refer to~\cite[Theorem 6.3.2]{dolgachev}
that
$b_i\, M_{\vartheta}\, b_j $ for distinct $i,j=1,\ldots,8$ define the 28 bitangents to $C$, arranged in Aronhold sets.
A recent result by Lehavi~\cite{lehavi} shows, that the set of bitangents uniquely determines $C$.
Moreover, there is a bijection between the 28 odd theta characteristics $\cL_{\vartheta_{ij}}$ and bitangents.
Any even theta characteristic different from $\cL_{\vartheta}$ can be represented by the divisor
class
\begin{equation} \label{othertheta}
\vartheta_{i,jkl} = \vartheta_{ij} + \vartheta_{ik} + \vartheta_{il} - K_C \mbox{ for distinct } i,j,k,l.\end{equation}

Next we define Aronhold bundles on $C$ following~\cite{pauly}. Given $\cJ \in\pic^1 C$ we define the 3--dimensional projective
space $\PP(\cJ):=\PP \ext^1(\omega_C \cJ^{-1},\cJ)=|\omega_C^2 \cJ^{-2}|^{\ast}$. A point in $\PP(\cJ)$
defines an isomorphism class of extensions
$$0\longrightarrow \cJ \longrightarrow \cK \longrightarrow \omega_C \cJ^{-1} \longrightarrow 0.$$
On $C$ pick the following data: 
\begin{itemize}
\item an even theta characteristic $\cL_{\vartheta}$,
\item a line bundle $\cJ_{\vartheta} \in\pic^1 C$ such that  $\cJ^2_{\vartheta}=\cL_{\vartheta}$,
\item a base point $b$ of the net of quadrics $\mathbb{M}_{\vartheta}$.
\end{itemize}
The stable (thus indecomposable) rank 2 bundle with canonical determinant $\cO_C(1)$ defined by the point $b\in \PP(\cJ_{\vartheta})$ is 
called the
\textit{Aronhold bundle} $\cK_{b,\vartheta}$. Up to 2--torsion points of $JC$, the bundles $\cK_{b,\vartheta}$ are in 
1-to-1 correspondence with the 288 unordered Aronhold sets. 

We mention a useful characterisation of Aronhold bundles: Let $\cK$ be a stable noneffective rank 2 vector bundle with
canonical determinant. By~\cite{lan} the set of line subbundles of maximal degree has cardinality 8 
$$\{\cJ \in\pic^1 C\ :\ \cJ \hookrightarrow \cK,\mbox{i.e. } h^0(\cJ^{-1} \otimes \cK)>0 \}=\{\cJ_1,\ldots,\cJ_8\}.$$
Since $\cK\in \PP(\cJ_i)$ for $i=1,\ldots,8$, there exist 28 effective divisors $D_{ij}$ of degree 2
satisfying
\begin{equation}\label{pauly28}
\cJ_i \otimes \cJ_j=\cO_C(D_{ij})\ \mbox{ for distinct } i,j=1,\ldots,8.
\end{equation}
Conversely, $\cK$ is uniquely determined by eight line bundles with property~(\ref{pauly28}) by~\cite{choe}. 
Finally, $\cK$ is an Aronhold bundle if and
only if the 28 effective divisors $D_{ij}$ correspond to bitangents on $C$.\\ 

Ottaviani~\cite{ottaviani} gives a nice description of the Aronhold invariant as a pfaffian which we briefly recall. 
The \textit{Aronhold invariant} of plane cubics is a quartic equation of $\sigma_3(\PP^2,\cO_{\PP^2}(3))$, which deals with the
condition to express an equation of a plane cubic as the sum of three cubes.
Here $$\sigma_3(\PP^2,\cO_{\PP^2}(3))\ =\mbox{ Zariski closure } \{g_1^3+g_2^3+g_3^3\ ;\ g_i \mbox{ linear forms} \} $$
denotes the 3-secant of the Veronese variety. 
Explicitly, the Aronhild invariant evaluated in 
\begin{eqnarray*}
w_{000}x^3+w_{111}y^3+w_{222}z^3+6w_{012}x y z+\\
3 w_{001}x^2 y+3w_{002}x^2 z+3w_{011}x y^2+3w_{022}x z^2+3w_{112}y^2 z
+3w_{122}y z^2
\end{eqnarray*}
equals $\pf Ar$ for
\begin{equation}
\label{aronhold}
Ar= \left[ \begin{array}{cccccccc}
0 & w_{222} & -w_{122} & 0 & w_{112} & 0 & w_{022} & -w_{012} \\
  & 0 & w_{022} & w_{122} & -w_{012} & -w_{022} & 0 & w_{002}\\
  &   & 0 & -w_{112} & 0 & w_{012} & -w_{002} & 0 \\
  &   &   & 0 & -w_{111} & 0 & -w_{012} & w_{011} \\
  &   &   &   & 0  & -w_{011} & w_{001} & 0 \\
  &   &   &   &    & 0 & w_{002} & -w_{001} \\
  &   &   &   &    &   & 0 & w_{000} \\
  &   &   &   &    &   &   &  0
  \end{array}\right].
  \end{equation}
On the other hand recall the  
\textit{Scorza map} between plane quartics~\cite{dolgachevkanev}:
$$F \mapsto \mbox{the Clebsch covariant quartic } S(F)$$
which is by definition
$$F \mapsto \mbox{polar cubic }P_{x}(F) \mbox{ at } x\in \PP^2 
\mapsto \mbox{Aronhold invariant}\! \left(P_{x}(F)\right).$$
Note that in this notation the coefficients $w_{ijk}$ of the cubic $P_{(x_0,x_1,x_2)}(F)$ are linear in $x_0,x_1,x_2$.
By~\cite[Section 7]{dolgachevkanev} the curve $S(F)$ carries a unique even theta characteristic $\vartheta$, more precisely,
the Scorza map
$$F \mapsto (S(F),\vartheta)$$
is an injective birational isomorphism and the natural projection to the first component is an unramified covering of degree 36.

\begin{proposition}\label{propott}
From the Aronhold pfaffian representation of $S(F)$ it is possible to recover the unique theta characteristic on $S(F)$.
\end{proposition}
\begin{proof}
Our main tool will be the Scorza correspondence relating points on $S(F)$. On
a nonsingular projective curve $X$ of genus $g > 0$ and a non--effective theta characteristic $\vartheta$ we introduce
\textit{ the Scorza correspondence }
$$R_{\vartheta}:=\{ (x,y)\in X\times X\, :\ h^0(\vartheta+x-y)>0  \}.$$
Dolgachev~\cite{dolgachev} proved that for $X$ non--hyperelliptic, $R_{\vartheta}$ are the only symmetric correspondences
of type $(g, g)$ without united points and some valence.
Therefore, on $S(F)$ there are 36 symmetric correspondences
of type $(3, 3)$.

We are able to explicitly determine which points on $S(F)$ are $R_{\vartheta}$ related in two essentially different ways:
\begin{itemize}
\item[(i)] from a symmetric determinantal representation of $S(F)$,
\item[(ii)] from an Aronhold pfaffian representation of $S(F)$.
\end{itemize}
The two ways which induce the same Scorza correspondence on $S(F)$ will relate the Aronhold pfaffian 
representation with the unique theta characteristic.

In step (i) we will use $M_{\vartheta}$,
the symmetric determinantal representation of $S(F)$ from Proposition~\ref{thetasymmdet}.
By definition $\lambda,\mu\in S(F)$ are NOT $R_{\vartheta}$ related if and only if $h^0(\vartheta+\lambda - \mu)=0$. This means that
$\cO_{S(F)}(\vartheta+\lambda-\mu) \oplus \cO_{S(F)}(\vartheta+\mu-\lambda)$ has canonical determinant and no sections.
Therefore, after tensoring by  $\cO_{S(F)}(1)$, it
equals the cokernel of another decomposable pfaffian representation of $S(F)$.
By Corollary~\ref{cordec} it
is obtained from 
$$(\cL_{\vartheta}\oplus \cL_{\vartheta})\otimes \cO_{S(F)}(1)\cong \coker M_{\vartheta} \oplus \coker M_{\vartheta}$$
by the Type I elementary transformation  based on the admissible vectors
$$\left[ \begin{array}{c}
v\\
0
\end{array}\right],\ 
\left[ \begin{array}{c}
0\\
u
\end{array}\right], \mbox{ where } v\in \coker M_{\vartheta}(\lambda),\, u\in \coker M_{\vartheta}(\mu).$$ 
The definition of admissible vectors~(\ref{vin6.7}) thus proves that $\lambda\, R_{\vartheta}\, \mu$ 
if and only if
$$v^t\, M_{\vartheta}(x)\, u\equiv 0 \mbox{ for all }
v\in \coker M_{\vartheta}(\lambda),\, u\in \coker M_{\vartheta}(\mu),\, x\in\PP^2.$$

Step (ii): Given an Aronhold pfaffian representation~(\ref{aronhold}), 
we can retrieve $F(x,y,z)$ from $w_{ijk}(x_0,x_1,x_2)$ by integrating 
$$P_{(x_0,x_1,x_2)}(F)=x_0 \frac{\partial F}{\partial x}(x,y,z)+x_1 \frac{\partial F}{\partial y}(x,y,z)+x_2 \frac{\partial F}{\partial z}(x,y,z).$$
The definition of the Aronhold invariant implies that, for any $\lambda\in S(F)$ there exist 
linear forms $g_1,g_2,g_3$ such that $P_{\lambda}(F)=g_1^3+g_2^3+g_3^3$.
This defines another $(3,3)$ correspondence without united points on $S(F)$, which must by the above equal to some $R_{\vartheta}$: 
$\lambda,\mu \in S(F)$ are \textit{related} if the second polar $P_{\lambda,\mu}(F)= g_i^2$ for some $i=1,2,3$. 
Obviously $\mu$ equals one of the vertices of the polar triangle spanned by the lines $g_1,g_2,g_3$. 
In~\cite[Theorem 7.8]{dolgachevkanev} Dolgachev and Kanev gave a beautiful construction 
of $F$ from the polar triangles in $S(F)$ and thus reconstructed $F$ from $(S(F),\vartheta)$.
\end{proof}

\begin{remark} {\rm Pauly's construction~\cite[\S 4.2]{pauly} assigns to every stable noneffective $\cK$ with
canonical determinant a net of quadrics whose bitangents correspond to $D_{ij}$ in~(\ref{pauly28}). This gives another proof of 
Proposition~\ref{propott} since the Aronhold bundle $\cK_{b,\vartheta}$ induces exactly the net of quadrics $\mathbb{M}_{\vartheta}$. 
We are however not able to implement this construction 
explicitly.}
\end{remark}

\begin{corollary} {\rm Denote by $R_{\vartheta}(\lambda)$ the polar triangle to $\lambda\in S(F)$. By~\cite{dolgachevkanev}, 
$R_{\vartheta}(\lambda)-\lambda$ equals the divisor class of $\cL_{\vartheta}$ and is thus independent of $\lambda$.
Then all the symmetric determinantal representations of $S(F)$ can be obtained from
$M_{\vartheta}$ by a sequence of three Type I elementary transformations, by applying the second part of Proof~\ref{cortwoparts} on
the divisor
$\vartheta_{i,jkl}-\vartheta=R_{\vartheta_{i,jkl}}(\lambda)-\lambda-R_{\vartheta}(\lambda)+\lambda=R_{\vartheta_{i,jkl}}(\lambda)-R_{\vartheta}(\lambda)$.
}
\end{corollary}

\begin{example} {\rm 
The Scorza map sends $F=x^4 + x^3 y - y^4 - y z^3 + 107^{1/3} x y^2 z$ to 
$$\begin{array}{ll}
S(F)=\pf [Ar]= &
27 x_0^3 x_1 - 432 x_0 x_1^3 - x_1^4 - 72\, 107^{1/3} x_0^2 x_1 x_2 - \\
  &  9\, 107^{1/3} x_0 x_1^2 x_2 + 81\, 107^{-1/3} x_0^2 x_2^2 - 108 x_0 x_2^3 - 27 x_1 x_2^3
 \end{array}$$
for $Ar$ defined in (\ref{aronhold}) with
$$\begin{array}{lllll}
w_{000}= 4 x_0 + x_1, & w_{001}= x_0,      & w_{011}= 1/3 x_2, & w_{111}= -4 x_1 , &  w_{002}= 0, \\
w_{012} = 1/3 x_1,    & w_{112} = 1/3 x_0, & w_{022} = 0,      & w_{122} =  - x_2 ,& w_{222} = -x_1.
\end{array}$$
Following the proof of Proposition~\ref{propott} we will compute the unique theta characteristic on $S(F)$. We calculate in
\texttt{Wolfram Mathematica} 
to precision $10^{-10}$.
For $\lambda = (1, 0, \frac{3}{4} 107^{-1/3} )\in S(F)$
we get 
$P_{\lambda}(F)=g_1^3+g_2^3+g_3^3$ for
$$\begin{array}{l}
g_1= (4 x + y) (-0.198` - 0.344`i),\\
g_2=(0.002` - 2.089` i) y + (-0.181` + 0.104` i) z ,\\
g_3= (0.002` + 2.089` i) y + (-0.181` - 0.104`  i) z, 
\end{array}$$ 
which is explicitly obtained from the equality $\det\, $Hess$\, (P_{\lambda}(F))=g_1 g_2 g_3$. 
The intersections 
$$\begin{array}{l}
\mu^1=g_2\cap g_3= (1,0, 0),\\   
\mu^2=g_1\cap g_3= (1, -4, -20.034` + 34.609` i),\\ 
\mu^3=g_1\cap g_2= (1, -4, -20.034` - 34.609` i)
\end{array}$$
determine the polar triangle $R(\lambda)$ of $\lambda$.
This proves that $\lambda$ is in relation with
$\mu^1,\ \mu^2$ and $\mu^3$ on $S(F)$.

On the other hand it is easy to compute all the 36 symmetric determinantal representations of $S(F)$. For example, for
$$\begin{array}{ll}
M_{\vartheta}= & x_1 \Id_4 - x_2 \mbox{ Diag}\, [\substack{0, -3, 3 (-1)^{1/3}, -3 (-1)^{2/3}}]+ \\
    &  x_0 
\left[\begin{array}{cccc}
\substack{4} & \substack{-24.296`} & \substack{23.685`+  0.336` i} & \substack{-23.685` + 0.336`i} \\
 &\substack{ \frac{428}{3}- 107^{1/3}} & \substack{-141.449` + 2.004` i} & \substack{141.449`+ 2.004` i}  \\
 &  & \substack{\frac{428}{3} - 107^{1/3}  (-1)^{2/3} }& \substack{-145.099`} \\
 &  &  & \substack{\frac{428}{3} + 107^{1/3}  (-1)^{1/3}}
\end{array}\right]
\end{array}$$
we have 
$$ \begin{array}{l}
v_{\lambda}=[\substack{-0.006 - 0.009 i, -0.335 - 
  0.482 i, -0.571 + 
  0.04 i, -0.236 + 0.521 i }]^t \in \coker M_{\vartheta}(\lambda),\\
  v_{\mu^1}=[\substack{ 0, -0.543 - 0.164 i, -0.419 + 0.404 i, 0.124+ 0.569 i}]^t \in \coker M_{\vartheta}(\mu^1),\\
  v_{\mu^2}=[\substack{ 0.602- 0.73 i, -0.186 -   0.025 i, -0.124 +  0.147 i, 0.062+ 0.172 i }]^t \in \coker M_{\vartheta}(\mu^2),\\
  v_{\mu^3}=[\substack{  0.613+ 0.72 i, -0.185 + 0.028 i, -0.059 + 0.173 i, 0.127+ 0.145 i }]^t \in \coker M_{\vartheta}(\mu^3).
  \end{array}$$
We check that
$$v_{\lambda}^t \, M_{\vartheta}(x_0,x_1,x_2)\, v_{\mu^i}=0\ \mbox{ for any }\ (x_0,x_1,x_2)\in\PP^2,\ i=1,2,3.$$
This proves that the corresponding $\cL_{\vartheta}$ is the unique theta characteristic on $S(F)$ that we were looking for.}
\end{example}

This is a counterexample to Ottaviani's conjecture~\cite[Remark 2.3]{ottaviani} that there exists an unique even characteristic $\cL_{\vartheta}$ 
on $S(F)$
for which $$h^0(\coker Ar \otimes \cO_{S(F)}(-1) \otimes \cL^{-1}_{\vartheta})>0.$$
Indeed, if $\coker Ar \otimes \cO_{S(F)}(-1) \otimes \cL^{-1}_{\vartheta}$ is effective, it is not stable since it has 
trivial determinant. Then $\coker Ar$ is also not stable and thus 
isomorphic to a direct sum of two line bundles. This is in contradiction with $Ar$ being an indecomposable representation of $S(F)$.


\end{document}